\documentclass[final]{siamltex}
\usepackage{amssymb}
\usepackage{amsmath}

\def\norm#1{\|#1\|}

\def\ZZ{\mathbb{Z}}
\def\RR{\mathbb{R}}
\def\ZZ{\mathbb{Z}}
\def\CC{\mathbb{C}}
\def\TT{\mathbb{T}}
\def\Cst{{\mathbb C}}

\def\Nst{{\mathbb N}}
\def\Qst{{\mathbb Q}}
\def\Rst{{\mathbb R}}

\def\Zst{{\mathbb Z}}
\def\Wil{{\mathcal{W}}}

\def\Hom{{\rm Hom}}

\def\fourier{{\mathcal F}}
\def\chirp{{\mathcal N}}
\def\dilate{{\mathcal D}}
\def\calA{{\mathcal A}}
\def\calJ{{\mathcal J}}

\def\theta{\vartheta}
\def\phi{\varphi}
\def\Lsp{{L}}
\def\lsp{{\ell}}

\def\LtR{{\Lsp^2(\Rst)}}
\def\ltZ{{\lsp^2(\Zst)}}

\def\vol{{\operatorname{vol}}}

\def\tc{{\tilde c}}
\def\td{{\tilde d}}

\def\tN{{N'}}
\def\Not{{\frac{N}{2}}}
\def\tNot{{\frac{\tN}{2}}}
\def\gab{{\cal G}}

\def\Lt{{\frac{L}{2}}}

\def\symp{\operatorname{Sp}(2,\Rst)}

\newcommand{\ip}[2]{\left\langle#1,#2\right\rangle}

\title{Wilson bases for general time-frequency lattices}

\author{Gitta Kutyniok\thanks{Institute of Mathematics, University of 
Paderborn, 33095 Paderborn, Germany; gittak@upb.de. G.K. 
acknowledges support from Forschungspreis 2003 der Uni\-ver\-si\-t\"at
Paderborn.} \and 
Thomas Strohmer\thanks{Department of Mathematics, University of
California, Davis, CA 95616-8633, USA; strohmer@math.ucdavis.edu. T.S.\
acknowledges support from NSF DMS grant 0208568.}}

\begin{document}

\maketitle

\begin{abstract}
Motivated by a recent generalization of the Balian-Low theorem 
and by new research in wireless communications
we analyze the construction of Wilson bases for general time-frequency 
lattices. We show that orthonormal Wilson bases for $\LtR$ can be 
constructed for any time-frequency lattice whose volume is $\tfrac12$. 
We then focus on the spaces $\ell^2(\ZZ)$ and $\CC^L$ which are the
preferred settings for numerical and practical purposes. We demonstrate
that with a properly adapted definition of Wilson bases the construction
of orthonormal Wilson bases for general time-frequency lattices also holds true
in these discrete settings. In our analysis we make use of certain metaplectic
transforms.
Finally we discuss some practical consequences of our theoretical findings.
\end{abstract}



\pagestyle{myheadings}
\thispagestyle{plain}
\markboth{G. KUTYNIOK AND T. STROHMER}{WILSON BASES FOR
GENERAL TIME-FREQUENCY LATTICES}

\section{Introduction} \label{intro}

Gabor systems have become a popular tool, both in theory and in applications,
e.g., see~\cite{FS98,Gro01,FS03}. However one drawback is that due to the 
{\em Balian--Low theorem} it is impossible to construct (orthogonal) Gabor 
bases for $\LtR$ with good time-frequency localization~\cite{Gro01}. 
In~\cite{Wil87,DJJ91} 
it has been shown that a modification of Gabor bases, so-called 
{\em Wilson bases} provide a means to circumvent the Balian--Low theorem. 
Indeed, there exist orthogonal Wilson bases for $\LtR$ whose basis functions 
have exponential decay in time and frequency. These Wilson bases can be 
constructed from certain tight Gabor frames with redundancy 2. 

Gabor frames are usually associated with rectangular time-frequency 
lattices, but they can also be defined for general non-separable lattices,
see e.g.~\cite{FCS95,Gro98,FK98,Kut00}. Recently it has been shown that
such a generalization of Gabor frames to general time-frequency lattices 
does not enable us to overcome the Balian--Low theorem~\cite{GHH02,BCM03}. 
This leads naturally to the question if it is possible to extend the
construction of Wilson bases to general time-frequency lattices. 

Another motivation for the research presented in this paper has its origin
in wireless communication.
Orthogonal frequency division multiplexing (OFDM) is a wireless
transmission technology employing a set of transmission functions which
is usually associated with a rectangular time-frequency
lattice~\cite{FAB95}. The connection to Gabor theory is given by the fact
that the collection of transmission pulses in OFDM can be interpreted
as a Gabor system, see~\cite{KM98,Str01}. 
The density of the associated rectangular time-frequency lattice can be 
seen as a measure of the spectral efficiency in terms of number of bits 
transmitted per Hertz per second. The necessary condition of linear 
independence of the transmission functions implies that we are dealing with 
either an undersampled or a critically sampled Gabor system.

For wireless channels that are time-dispersive (due to multipath) and
frequency-dispersive (due to the Doppler effect) good time-frequency 
localization of the transmission pulses is essential to mitigate the
interferences caused by the dispersion of the channel~\cite{KM98,Str01}. The ideal set of
transmission pulses should therefore possess  (i) good time-frequency 
localization and (ii) maximize the spectral efficiency, i.e., the transmission
functions should correspond to a (critically sampled) Gabor basis
for $\LtR$. As we know, the Balian--Low theorem prohibits these conditions 
to be fulfilled simultaneously. 

Recently it has been shown that in case of
time-frequency dispersive channels the performance of OFDM can be improved
when using general time-frequency lattices, in particular hexagonal-type
lattices~\cite{SB03}. In a nutshell, lattices that are adapted to the shape 
of the Wigner distribution of the transmission pulses allow for 
a better ``packing'' of the time-frequency plane, which in turn can be used
to either achieve higher data rates or to improve interference robustness
of the associated so-called Lattice-OFDM system.

One variation of OFDM (for rectangular lattices) is called Offset-QAM (OQAM)
OFDM, it corresponds to using a Wilson basis as set of transmission 
functions~\cite{Bol02}. 
OQAM-OFDM achieves maximal spectral efficiency and allows for transmission
functions with good time-frequency localization.
As in the case of standard OFDM it would be potentially useful for 
time-frequency dispersive channels to extend OQAM-OFDM to general 
time-frequency lattices in order to improve the robustness of  
OQAM-OFDM against interference even further. Thus we again arrive at the 
problem of constructing Wilson bases for general non-separable 
time-frequency lattices. 

Yet another motivation comes from filter bank theory, more precisely
cosine-modulated filter banks~\cite{BH98}. We know that discrete-time Wilson 
bases correspond to a special class of cosine-modulated filter banks 
(see~\cite{BH98}).
In light of the improvements gained by using general time-frequency lattices 
in OFDM \cite{SB03}, it would be interesting to analyze if the construction 
of cosine-modulated filter banks can be extended to general time-frequency 
lattices. A positive answer to this question might lead to a more efficient 
encoding of signals and images. 

Since our goal to construct Wilson bases for general time-frequency lattices 
is in part motivated by applied problems and since any numerical 
implementation of Wilson bases is based on a discrete model, our analysis will
not only concern $\LtR$ but also comprise the spaces $\ltZ$ and $\Cst^L$.
Furthermore, $\ltZ$ is the appropriate setting when Wilson bases are 
utilized as filter banks, since in this case one deals with
sampled, thus discrete-time, signals.

\subsection{Notation} \label{ss:not}

We assume that the reader is familiar with the theory of Gabor frames and
refer to~\cite{Gro01} for background and details.

A lattice $\Lambda$ in $\Rst^d$ is a discrete subgroup with compact
quotient, i.e., there exists a matrix $A\in GL(d,\Rst)$ such that
$\Lambda = A \Zst^d$. The matrix $A$ is called the (non-unique)
generator matrix for $\Lambda$. The volume of $\Lambda$ is $\vol (\Lambda) = 
|\det(A)|$. Two lattices, which play a crucial role in OFDM design (see
\cite{SB03}), are the rectangular lattice $\Lambda_R$ and the hexagonal lattice
$\Lambda_H$. A generator matrix for $\Lambda_R$ is given by
\[ A_R = \begin{bmatrix}
T & 0 \\
0 & F 
\end{bmatrix}\]
and a generator matrix for $\Lambda_H$ is given by
\[ A_H = \begin{bmatrix}
\frac{\sqrt{2}}{\sqrt[4]{3}} T & \frac{\sqrt{2}}{2\sqrt[4]{3}} T \\
0 & \frac{\sqrt[4]{3}}{\sqrt{2}} F 
\end{bmatrix},\]
where $T,F > 0$. An easy calculation shows that both lattices $\Lambda_R$
and $\Lambda_H$ have the same volume $TF$.
A normal form for matrices, which we will use in the following, is the so-called
Hermite normal form \cite{Her1850}. 
We say that a matrix $A = 
\begin{bmatrix}
a & b \\
c & d 
\end{bmatrix}$
is in {\em Hermite normal form}, if $c=0$, $a,d > 0$, and $0 \le b < a$. For example
both matrices $A_R$ and $A_H$ are in Hermite normal form.

For $(x,y) \in \mathbb{R}^2$ and $g \in L^2(\mathbb{R})$, let $g_{x,y}$ be
defined by
\[g_{x,y}(t) = e^{2 \pi i y t} g(t-x).\]
We denote by $\gab(g,\Lambda)$ the system of functions given by
\begin{equation*}
g_{\lambda,\mu}(t)=g(t-\lambda) e^{2\pi i t\mu},
\qquad (\lambda,\mu) \in \Lambda.
\end{equation*}
As usual, the {\em redundancy} of $\gab(g,\Lambda)$ is given by
$\frac{1}{\vol(\Lambda)}$.

As in \cite{Gro01}, we define the {\em Schr\"odinger representation}
$\rho : \mathbb{H} \to \mathcal{U}(L^2(\mathbb{R}))$  by
\begin{equation*}
\rho(x,y,z)g(t) = e^{2 \pi i z}e^{-\pi i xy} e^{2 \pi i y t} g(t-x).
\end{equation*}
Note that
\begin{equation}\label{xy} g_{x,y} = e^{\pi i x y} \rho(x,y,1)g.\end{equation}
Furthermore, we use the following notation from \cite{Gro01} (with slight
changes):
\[\mathcal{J} =  \begin{bmatrix} 0 & -1 \\ 1 & 0
\end{bmatrix}, \quad
\mathcal{B}_b =  \begin{bmatrix} b & 0 \\ 0 & \frac{1}{b}
\end{bmatrix}, \quad  
\mathcal{C}_c =  \begin{bmatrix} 1 & 0 \\ c & 1
\end{bmatrix}.\]
$\mathcal{F}$ and $^\wedge$ denote the Fourier transform, 
$\mathcal{F}^{-1}$ and $^{\vee}$ denote the inverse Fourier transform.
The dilation is given by $\mathcal{D}_b f(t) =
|b|^{\frac{1}{2}} f(bt)$ and the ``chirp'' operator is defined via
$\mathcal{N}_c f (t) = e^{-\pi i c t^2} f(t)$,
where $b,c \in \mathbb{R}$ and $f \in L^2(\mathbb{R})$.

Metaplectic transforms will turn out to be a very useful tool in our analysis.
For the study of the discrete and finite case, we need a result on
metaplectic transforms from
\cite{Kai99}. Since this thesis is not easy to access, we present  
the result here in the slightly weaker version we will use, together with the necessary
definitions and notations. The result will be stated in the  situation
of a general locally compact abelian group $G$ with dual group $\widehat{G}$,
group multiplication denoted by $+$,
and action of $\widehat{G}$ on $G$ denoted by $\ip{x}{\chi}$ for
$x\in G$ and $\chi\in\widehat{G}$. 
The cases we are interested in later on are $G=\ZZ$ and $G=\ZZ_L$. 
We will usually write a metaplectic transform $\sigma$ in the matrix notation 
$\sigma = \left[ \begin{array}{cc} \alpha & \beta\\ \gamma & \delta \end{array}\right] 
\in \Hom(G \times \widehat{G})$, which means that $\alpha \in \Hom(G)$,
$\beta \in \Hom(\widehat{G},G)$, $\gamma \in \Hom(G,\widehat{G})$, and 
$\delta \in \Hom(\widehat{G})$.
Then the {\em adjoint} $\sigma^*$ is defined by $\ip{(x,\chi)}{\sigma^*(y,\pi)}
= \ip{\sigma(x,\chi)}{(y,\pi)}$ for all $(x,\chi),(y,\pi) \in G \times \widehat{G}$.
Let $\eta$ be defined by
$\eta = \left[ \begin{array}{cc} 0 & -I_{\widehat{G}} \\I_G & 0 \end{array}\right]
\in \Hom(G \times \widehat{G},\widehat{G}
\times G)$,
where the above definition concerning the matrix notation has to be adapted in 
an obvious way, and 
where $I_G$ and $I_{\widehat{G}}$ denote the identity on $G$ and $\widehat{G}$,
respectively.
Then $\sigma$ is called {\em symplectic}, if $\sigma^* \eta \sigma = \eta$.
If $\zeta \in \Hom(G \times \widehat{G},\widehat{G}
\times G)$, then $\psi$ is a {\em second degree character} of $G \times \widehat{G}$
associated to $\zeta$, if
$\psi((x,\chi)+(y,\pi)) = \psi(x,\chi)\psi(y,\pi)
\ip{(x,\chi)}{\zeta(y,\pi)}$ for all $(x,\chi),(y,\pi) \in G \times \widehat{G}$.
Moreover, for $(x,\chi) \in G \times \widehat{G}$ and $g \in L^2(G)$, let
$g_{x,\chi}$ be defined by
\[g_{x,\chi}(t) = \chi(t)g(t-x),\]
which generalizes the previous definition.

We can now state the version of \cite[Theorem~1.1.28]{Kai99}, which we will 
employ in Sections \ref{s:disc} and \ref{s:fin}.

\begin{theorem}\label{Kai_Theorem}
Let $\sigma := \left[ \begin{array}{cc} \alpha & \beta\\
\gamma & \delta \end{array}\right] \in \Hom(G \times \widehat{G})$ be symplectic
and let $\psi_\zeta$ be a second degree character of $G \times \widehat{G}$
associated to 
\[\zeta := \sigma^* \left[ \begin{array}{cc} 0 & 0\\I_G & 0 \end{array}\right] \sigma 
- \left[ \begin{array}{cc} 0 & 0\\I_G & 0 \end{array}\right] \in \Hom(G \times \widehat{G},\widehat{G}
\times G).\]
If $U$ is defined by
\[Uf(t):= \int_{\widehat{G}} f(\alpha t + \beta \omega) \psi_\zeta^{-1}(t,\omega) d\omega,\]
then we have
\[ (Uf)_{x,\chi}(t) =  \psi_\zeta^{-1}(x,\chi) U f_{\sigma(x,\chi)}(t), \quad (x,\chi) \in 
G \times \widehat{G}.\]
\end{theorem}


\section{Wilson bases for general lattices -- the continuous case} 
\label{s:cont}

We first show that all lattices in $\RR^2$ which are important for applications,
such as the rectangular lattice, the hexagonal lattice, and lattices whose
generator matrix has rational entries, possess a uniquely determined
matrix in Hermite normal form, which we will call the {\em canonical generator
matrix}. In particular, we characterize exactly those lattices, which possess a
generator matrix in Hermite normal form. 

We then show that it is possible 
to construct orthonormal Wilson bases for time-frequency lattices 
$\Lambda$ with $\vol (\Lambda) = \frac{1}{2}$, which possess a 
generator matrix in Hermite normal form.
In principle Wilson systems can be defined for lattices $\Lambda$ with
$\vol(\Lambda)\neq \frac{1}{2}$, however so far all known constructions of 
{\em orthogonal Wilson bases} for $\LtR$ are strictly tied to lattices 
with volume $\frac{1}{2}$. In light of this fact throughout this paper a 
Wilson system will always be associated with a time-frequency lattice of 
volume $\frac{1}{2}$.

\begin{lemma} \label{canonicalgeneratormatrix}
Let $\Lambda$ be a lattice in $\RR^2$. Then the following conditions are equivalent.
\begin{enumerate}
\item[\rm (i)] $P_2(\Lambda)$ is discrete, where
$P_2 : \RR^2 \to \RR$, $(x,y) \mapsto y$.
\item[\rm (ii)] There exists a generator matrix $A$ for $\Lambda$ which is in Hermite normal form.
\end{enumerate}
If one of these conditions is satisfied, the matrix $A$ is uniquely determined.
\end{lemma}

\begin{proof}
Let 
\[ A' = 
\begin{bmatrix}
a' & b' \\
c' & d' 
\end{bmatrix}\]
be an arbitrary generator matrix for $\Lambda$. 

First we prove that~(ii) implies~(i). By~(ii), there exists a matrix
\[ A = 
\begin{bmatrix}
a & b \\
0 & d 
\end{bmatrix},\]
which is in Hermite normal form and which satisfies $A \ZZ^2 = \Lambda$. Thus
$P_2(\Lambda) = d \ZZ$, which yields~(i). 

Next we show~(i) $\Rightarrow$~(ii). For this, we construct a matrix
\[ A = 
\begin{bmatrix}
a & b \\
0 & d 
\end{bmatrix},\]
which satisfies the claimed properties. Without loss of generality we assume that $d' \neq 0$
(if $d' = 0$ and $c' \neq 0$ we could change the columns of $A'$).
We begin with the following observation.
Assume that $\frac{c'}{d'}$ is not rational. Since $P_2(\Lambda)$ is a non-trivial discrete,
additive subgroup of $\RR$, hence a lattice, there exists $s \in \RR \backslash \{0\}$ 
such that $P_2(\Lambda)=s\ZZ$. Thus $c'=sm$ and $d'=sn$ for some $m,n \in \ZZ$, a contradiction.
This implies that the quotient $\frac{c'}{d'}$ is rational. Setting $\frac{c'}{d'} = \frac{k}{l}$, 
$k,l \in \ZZ$ with $\gcd(k,l)=1$
and factoring out $\frac{d'}{l}$, without loss of generality we can assume that
$A'$ is of the form
\[ A' = r
\begin{bmatrix}
a' & b' \\
c' & d' 
\end{bmatrix},\]
with $a',b',r \in \RR$ and $c',d' \in \ZZ$. Now we proceed as follows.
First we set $p := \gcd(c',d')>0$. 
We then obtain
\[ A'
\begin{bmatrix}
\pm \frac{d'}{p} \\
\mp \frac{c'}{p}
\end{bmatrix}
= \begin{bmatrix}
\pm \frac{r\det(A')}{p} \\
0
\end{bmatrix} \in \Lambda,\]
since $\frac{c'}{p}, \frac{d'}{p} \in \ZZ$. 
Hence we can define $a$ by  $a := \left|\frac{r\det(A')}{p}\right| > 0$.  

In a second step we compute $b$ and $d$. Since $\frac{c'}{p}$ and  $\frac{d'}{p}$ are
relative prime, there exist $m,n \in \ZZ$ such that $\frac{c'}{p}m + \frac{d'}{p}n=1$
(see \cite[Theorem~4.4]{Hua82}). Hence $c'm+d'n=p$ and we obtain
\[ A'
\begin{bmatrix}
m \\
n
\end{bmatrix}
= \begin{bmatrix}
r(a'm+b'n) \\
rp
\end{bmatrix} \in \Lambda.\]
Now let $k \in \ZZ$ be chosen in such a way that
$0 \le r(a'm+b'n)+ka   < a$ and define $b := ka + r(a'm+b'n)$ and $d := rp$. Without loss
of generality we can assume that $d > 0$, since otherwise we just take $-m$ and $-n$
instead of $m$ and $n$. Then
\[ \begin{bmatrix}
b \\
d
\end{bmatrix}
= \begin{bmatrix}
r(a'm+b'n) \\
rp
\end{bmatrix}
+ k\begin{bmatrix}
a \\
0
\end{bmatrix}
\in \Lambda\]
and $|ad| = \left|\frac{r\det(A')}{p}\right|rp = |r^2\det(A')|$. This proves that $A$ generates
$\Lambda$ and is in Hermite normal form.

At last we prove that the matrix in condition~(ii) is uniquely determined. For this, assume there 
exist $a,b,d,a',b',d' \in \RR$ with $a,a',d,d' > 0$, $0 \le b < a$ and
$0 \le b' < a'$ such that 
\begin{equation} \label{equal}
\Lambda = A \ZZ^2 = A' \ZZ^2,
\end{equation}
where
\[A = 
\begin{bmatrix}
a & b \\
0 & d 
\end{bmatrix}
\quad \mbox{and} \quad
A' = 
\begin{bmatrix}
a' & b' \\
0 & d' 
\end{bmatrix}.\]
By (\ref{equal}), there exist $m,n \in \ZZ$ with $a' = ma+nb$ and $0 = nd$, which
implies that $a'=ma$. Again by (\ref{equal}),  we can find $k,l \in \ZZ$
such that $b'=ka+lb$ and $d'=ld$. Since $\vol(\Lambda) = |ad| = |a'd'|$, we obtain
$|ml|=1$. Now  $a,a',d,d' > 0$ implies that $a = a'$, $d = d'$, and $l=1$. Finally,
applying this to $b'=ka+lb$ and using that $0 \le b < a$ and
$0 \le b' < a'$ yields $b=b'$. Thus we have shown $A=A'$, which completes the proof.
\end{proof}

In the following we restrict our attention to lattices, which possess a
generator matrix in Hermite normal form. All results in the situation $L^2(\RR)$ 
could be derived (in the same manner, but with much more technical effort) for general 
lattices, however with little or no practical benefit.

\begin{definition}
Let $\Lambda$ be a lattice in $\RR^2$, which satisfies the conditions of Lemma
\ref{canonicalgeneratormatrix}. Then the uniquely determined
generator matrix $A$ of Lemma \ref{canonicalgeneratormatrix} is called the
{\em canonical generator matrix} for $\Lambda$. 
\end{definition}

Now let $\Lambda$ be a lattice with $\vol(\Lambda)=\tfrac12$, which possesses a
generator matrix in Hermite normal form.
Using the definition of a canonical generator matrix, we
define a Wilson system associated with $\Lambda$ as follows.
\begin{definition}
If $\gab(g,\Lambda) \subseteq \LtR$ is a Gabor system of redundancy 2, and
\begin{equation*}
A = 
\begin{bmatrix}
a & b \\
0 & d 
\end{bmatrix}
\end{equation*}
is the canonical generator matrix for the lattice $\Lambda$, then
the associated {\em Wilson system} $\Wil(g,\Lambda,L^2(\RR)) =
\{\psi_{m,n}^\Lambda\}_{m\in\Zst,n\ge 0}$ consists of the functions 
\begin{eqnarray*}
\psi_{m,0}^\Lambda & = & g_{2ma,0},\hspace*{5.02cm}\mbox{if $n =0$}, \\
\psi_{m,n}^\Lambda & = & \tfrac{1}{\sqrt{2}}  e^{-\pi i bdn^2}
(g_{ma+nb,nd}+g_{ma-nb,-nd}),\quad
\mbox{if } n \neq 0, \: m+n \mbox{ even},\\
\psi_{m,n}^\Lambda & = & \tfrac{i}{\sqrt{2}} e^{-\pi i bdn^2}
(g_{ma+nb,nd}-g_{ma-nb,-nd}),\quad
\mbox{if } n \neq 0, \: m+n \mbox{ odd}.
\end{eqnarray*}
If the system $\Wil(g,\Lambda,L^2(\RR))$ is an orthonormal basis for $\LtR$ we call it
a Wilson (orthonormal) basis.
\end{definition}

We will see that this definition reduces to the usual definition of Wilson systems.
For this, we consider the rectangular lattice 
\[\Gamma = \{(\tfrac{m}{2},n)\}_{m,n \in \Zst}.\]
It is an easy calculation to show that the canonical generator matrix for $\Gamma$ is 
\[ A = \begin{bmatrix}
\frac12 & 0 \\
0 & 1 
\end{bmatrix}.\]
Thus, for each $g \in L^2(\RR)$, the Wilson system $\Wil(g,\Gamma,L^2(\RR))$
consists indeed of the functions
\begin{eqnarray*}
\psi_{m,0}^{\Gamma}& = & g_{m,0}, \hspace*{3cm} \text{if $n = 0$},  \\
\psi_{m,n}^{\Gamma} & = & \tfrac{1}{\sqrt{2}}(g_{m/2,n}+g_{m/2,-n}), \quad 
\text{if $m+n$ is even}, \\
\psi_{m,n}^{\Gamma} & = & \tfrac{i}{\sqrt{2}}(g_{m/2,n}-g_{m/2,-n}), \quad 
\text{if $m+n$ is odd}, 
\end{eqnarray*}
which coincides with the usual definition of Wilson systems,
cf.~for instance \cite[Definition~8.5.1]{Gro01}. Notice that we will fix the
notation $\Gamma$ for the remainder.

We will make use of the following well-known theorem about a Wilson system
for rectangular lattices to constitute an orthonormal basis 
(e.g.~cf. \cite[Theorem~4.1]{BHW95}).
\begin{theorem}
\label{th:BHW}
Suppose that $g \in  \LtR$ is such that
\begin{itemize}
\item[\rm (a)] $\hat{g}$ is real-valued and
\item[\rm (b)] $\{g_{m/2,n}\}_{m,n\in \Zst}$ is a tight Gabor frame 
for $\LtR$ with frame bound $2$.
\end{itemize}
Then the system $\Wil(g,\Gamma,L^2(\RR))$
is a Wilson orthonormal basis for $\LtR$.
\end{theorem}


We are now ready to extend the construction of Wilson bases for
time-frequency lattices which possess a generator matrix in Hermite 
normal form.

\begin{theorem}
\label{theo1}
Let $\Lambda$ be a lattice in $\Rst^2$ with $\vol (\Lambda) = \frac{1}{2}$
and canonical generator matrix 
\begin{equation*}
A = 
\begin{bmatrix}
a & b \\
0 & d 
\end{bmatrix}.
\end{equation*}
Define $U$ by
\begin{equation}
U : =  \dilate_{1/d} \circ
\fourier \circ \chirp_{-b/d} \circ \fourier^{-1}.
\label{Umatrix}
\end{equation}
Let $g \in L^2(\mathbb{R})$ be such that
\begin{itemize}
\item[\rm (i)] $\widehat{Ug}$ is real-valued,
\item[\rm (ii)] $\{g_{ma+nb,nd}\}_{m,n \in \mathbb{Z}}$ is a tight frame for 
$\LtR$ with frame bound $2$.
\end{itemize}
Then the system $\Wil(g,\Lambda,L^2(\RR))$
is a Wilson orthonormal basis for $L^2(\mathbb{R})$.
\end{theorem}

\begin{proof}
We will reduce our claim to Theorem~\ref{th:BHW} by using a 
metaplectic transform. We define $\mathcal{A}$ by
\begin{equation*}
\mathcal{A} =
\begin{bmatrix}
d & -b \\
0        & 2 a
\end{bmatrix}.
\end{equation*}
Then, for all $m,n \in \mathbb{Z}$, we have
\begin{equation}
\label{latticeconnetion}
\calA (ma+nb,nd) = \big(\tfrac{1}{2}m,n\big).
\end{equation}
Since $\calA \in \symp$ we can apply \cite[Theorem~4.51]{Fol89}  and write
\[\calA = 
 B_{d} (-\calJ) C_{b/d} \calJ.\]
By \cite[Example~9.4.1]{Gro01} we obtain
\begin{equation} \label{rhoequation}
\rho(x,y,1)g =  U^{-1}\rho(\mathcal{A}(x,y),1) (Ug),
\end{equation}
with $U$ defined in~\eqref{Umatrix}. Using (\ref{xy}), (\ref{latticeconnetion}), 
and (\ref{rhoequation}), we compute 
\begin{eqnarray*}
g_{ma+nb,nd}& = & e^{\pi i (ma+nb)nd} \rho(ma+nb,nd,1)g \\
 & = & e^{\pi i (ma+nb)nd} U^{-1}\rho(\calA(ma+nb,nd),1)(Ug)\\
 & = & e^{\pi i (ma+nb)nd} U^{-1}\rho(\tfrac{1}{2} m,n,1)(Ug). 
\end{eqnarray*}
Using (\ref{xy}) again, we obtain
\begin{equation} \label{relation}
g_{ma+nb,nd} = e^{\pi i (ma+nb)nd} e^{-\pi i \frac{1}{2} mn} U^{-1}(Ug)_{m/2,n}
= e^{\pi i bdn^2} U^{-1}(Ug)_{m/2,n}.
\end{equation}

Since multiplication by a phase factor and applying a unitary operator to a 
tight frame preserves tightness (and frame bounds), it follows that condition~(ii) 
is equivalent to $\{(Ug)_{m/2,n}\}_{m,n\in \Zst}$
being a tight frame with frame bound $2$. We need not deal with condition (i),
since this states already that $\widehat{Ug}$ is real-valued. Applying
Theorem~\ref{th:BHW} yields that the Wilson system  $\Wil(Ug,\Gamma,L^2(\RR))$
is an orthonormal basis. Using now the metaplectic transform, i.e., (\ref{relation}),
and the fact that $U$ is a unitary operator, finishes the proof.
\end{proof}

We will conclude this section by providing an example, which shows how
we can compute a Wilson basis with excellent time-frequency localization
for a special lattice, but this calculation can also be done for an arbitrary lattice.

Here we consider the hexagonal lattice $\Lambda_H$ with generator matrix
\[ A_H = \begin{bmatrix}
\frac{\sqrt{2}}{\sqrt[4]{3}} & \frac{\sqrt{2}}{2\sqrt[4]{3}} \\
0 & \frac{\sqrt[4]{3}}{\sqrt{2}}  
\end{bmatrix}.\]
This lattice was also used in~\cite{SB03}. Observe that since 
we have $0 \le \frac{\sqrt{2}}{2\sqrt[4]{3}} < \frac{\sqrt{2}}{\sqrt[4]{3}}$, 
the matrix $A_H$ is already the canonical generator matrix for $\Lambda_H$. 
First we define the function $h \in L^2(\RR)$ by
\begin{equation*}
h(x) = (2\nu)^{\frac{1}{4}}e^{-\nu\pi x^2}.
\end{equation*}
By \cite[Theorem~7.5.3]{Gro01}, the set $\{h_{m/2,n}\}_{m,n \in \ZZ}$ is a frame
for $L^2(\RR)$. Let $S$ denote its frame operator and consider the
function $\phi \in L^2(\RR)$ given by
\begin{equation}
\label{squareroot}
\phi = \sqrt{2} \mathcal{F} \circ S^{-\frac12}h.
\end{equation}
Using \cite[Theorem~4.6]{BHW95}, this function coincides with the function considered
in \cite[Section~4]{DJJ91}. There it was shown that $\phi$ satisfies conditions (a) 
and (b) of Theorem~\ref{th:BHW} and hence yields a Wilson basis in the sense
of Theorem~\ref{th:BHW}. Moreover, $\phi$ has exponential decay in time
and frequency. To obtain a generating function for a Wilson basis
with respect to $\Lambda_H$ in the sense of Theorem~\ref{theo1}, first 
observe that, by the proof of Theorem~\ref{theo1}, we only need to compute
the function $g=U^{-1}\phi$, where $U$ is defined in (\ref{Umatrix}). Then
$g$ automatically satisfies conditions (i) and (ii) of Theorem~\ref{theo1},
and hence the system $\Wil(g,A_H,L^2(\RR))$ is a Wilson orthonormal basis by
Theorem~\ref{theo1}. Thus we define $g \in L^2(\RR)$ by
\[ g = \mathcal{F} \circ \mathcal{N}_{\frac{1}{\sqrt{3}}} \circ \mathcal{F}^{-1}
\circ \mathcal{D}_{\frac{\sqrt[4]{3}}{\sqrt{2}}} \phi
=  \sqrt{2} \mathcal{F} \circ \mathcal{N}_{\frac{1}{\sqrt{3}}} \circ 
\mathcal{D}_{\frac{\sqrt{2}}{\sqrt[4]{3}}} \circ S^{-\frac12}h.\]
Let us mention that the function $g$ has exponential decay in time
and frequency. Thus we obtain a Wilson basis with respect to the lattice
$\Lambda$ with very good time-frequency localization. As already mentioned
above this procedure can be applied to an arbitrary lattice, hence we obtain a
Wilson basis with excellent time-frequency localization for any lattice.


\section{Wilson bases for general lattices -- the discrete case} 
\label{s:disc}

In this section we analyze the construction of Wilson bases for
general time-frequency lattices for functions defined on $\ltZ$. 
The reasons for considering the setting $\ltZ$ are on that the one hand 
several applications such as filter bank design in digital signal
processing deal directly with a discrete setting~\cite{BH98}, and on the other
hand, even those problems that arise in the ``continuous'' setting
of $\LtR$ require a discrete model for their numerical treatment.
Thus, with these practical aspects in mind, throughout this section 
we naturally consider only 
lattices whose generator matrices have rational entries, since any
implementation is intrinsically restricted to such ``rationally'' 
generated lattices. Another natural setting for numerical implementations is 
of course $\Cst^L$ (which can be identified with the space of
$L$-periodic sequences). We will analyze that case in the next section.

Before we proceed we define Gabor systems and Wilson systems on $\ltZ$ for 
general time-frequency lattices $\Lambda$ with $\vol(\Lambda)=\frac{1}{2}$.
First we prove that each lattice possesses a generator matrix of some particular
form.

\begin{lemma}
\label{lemmadisc}
Let $\Lambda$ be a lattice in $\ZZ \times \RR$ with generator matrix $A$ given by
\begin{equation*}
A=\begin{bmatrix} a & b \\ c & d \end{bmatrix},
\qquad \text{with $a,b\in \Zst$, $c,d \in \Qst$, and $\det(A)=\frac{1}{2}$},
\end{equation*}
and denote $c=\frac{p}{q}, d=\frac{p'}{q'}$ with $\gcd(p,q)=\gcd(p',q')=1$,
$p,q,p',q'\in \Zst$. 
Then $\Lambda$ possesses a uniquely determined generator matrix of the form
\begin{equation}
A' =
\begin{bmatrix}
\frac{N}{2} & b' \\
0  & \frac{1}{N}
\end{bmatrix},
\label{redA2}
\end{equation}
where $N=\frac{qq'}{\gcd(pq',p'q)}$ and $b' \in \Zst$, $0 \le b' < \frac{N}{2}$. 
\end{lemma}

\begin{proof}
We assume that $c\neq 0$, otherwise~\eqref{redA2} 
is automatically satisfied. We first show that $A$ can be written as
\begin{equation}
A =
\begin{bmatrix}
a & b \\
\frac{r}{N} & \frac{s}{N}
\end{bmatrix},
\label{redA3}
\end{equation}
with integers $r,s,N$, such that 
$\gcd(r,s)=1$. 

Let $c=\frac{p}{q}, d=\frac{p'}{q'}$ with $p,p',q,q' \in \Zst$, denote
$\tN:=qq', \tc: = pq', \td: =p'q$ and write $z = \gcd(\tc,\td)$. Since
 $a,b,\tc,\td \in \Zst$ and since
$$
\vol(\Lambda) = \frac{1}{2} \Rightarrow a\td - b\tc = \frac{\tN}{2}, 
$$
it follows that $\frac{\tN}{2} \in \Zst$. A necessary and sufficient
condition for the equation $a\td - b \tc = \frac{\tN}{2}$ to have an integer 
solution in $a$ and $b$ is that $\gcd(\tc,\td)|\frac{\tN}{2}$, see
\cite[Theorem~8.1]{Hua82}, hence $z|\tNot$. Denote
$z' :=\frac{N'}{2z} ,r : = \frac{\tc}{z}, s:=\frac{\td}{z}$. Then
$c= \frac{r}{2z'}, d=\frac{s}{2z'}$ with
$z' \in \Zst$ and $\gcd(r,s)=1$. By a proper choice of the signs
of $c$ and $d$ we can always assume that $z'$ is positive. By writing
$N:=2z' \in \Zst$ we see that $A$ can indeed be written as 
in~\eqref{redA3}. 

Now, assuming that $A$ is of the form~\eqref{redA3}, we compute
\begin{equation*}
\begin{bmatrix}
a & b \\
c & d
\end{bmatrix}
\begin{bmatrix}
Nd \\-Nc 
\end{bmatrix}
=
\begin{bmatrix}
N(ad-bc)\\0 
\end{bmatrix}
=
\begin{bmatrix}
\Not \\ 0
\end{bmatrix}.
\end{equation*}
Since $\gcd(r,s) =1$ there exist integers $m,n$ with
$mr +ns = 1$. For such a pair $(m,n)$ we denote $b'=am+bn$ and obtain
\begin{equation*}
\begin{bmatrix}
a & b \\
c & d
\end{bmatrix}
\begin{bmatrix}
m \\n  
\end{bmatrix}
=
\begin{bmatrix}
am+bn \\ \frac{r}{N} m +\frac{s}{N} n
\end{bmatrix}
=
\begin{bmatrix}
b' \\ \frac{1}{N}
\end{bmatrix}.
\end{equation*}
If $b' < 0$ or $b' \ge \frac{N}{2}$, we substitute the vector obtained by
\[ \begin{bmatrix}
b' \\ \frac{1}{N}
\end{bmatrix}
+k
\begin{bmatrix}
\Not \\ 0
\end{bmatrix}
= \begin{bmatrix}
b'+k\Not \\ \frac{1}{N}
\end{bmatrix},\]
where $k \in \ZZ$ is chosen in such a way that $0 \le b'+k\Not < \frac{N}{2}$.
Consequently the matrix 
\begin{equation}
\begin{bmatrix}
\Not & b' \\
0 & \frac{1}{N} 
\end{bmatrix}
\label{redA4}
\end{equation}
generates the lattice $\Lambda$. Finally, since $\Not$ and $b'$ are 
integers and $0 \le b' < \frac{N}{2}$, the generator matrix in~\eqref{redA4} 
is indeed of the form~\eqref{redA2}. 

The fact that this is a unique representation follows immediately from
the condition $0 \le b' < \frac{N}{2}$.
\end{proof}

\begin{definition}
Let $\Lambda$ be a lattice in $\ZZ \times \RR$. Then the uniquely determined
matrix $A'$ of Lemma \ref{lemmadisc} is called the
{\em canonical generator matrix} for $\Lambda$. 
\end{definition}

In the following we will regard such a lattice as a lattice in 
$\ZZ \times \TT$ by considering 
$\Lambda = \{\frac{N}{2}m+bn,e^{2 \pi i \frac{n}{N}}\}_{m\in\ZZ,
n=0,\ldots,N-1}$. This is a very natural approach, since, for all $k \in \ZZ$, we have
\[ (\tfrac{N}{2}m+bn,\tfrac{1}{N}n+k)
= (\tfrac{N}{2}(m-2bk)+b(n+kN),\tfrac{1}{N}(n+kN)),\]
hence the lattice $A' \ZZ^2$ is invariant under adding integers to the 
second component.  Moreover, it is sufficient to restrict to the
index set $\ZZ \times \{0,\ldots,N-1\}$, since, for all $0\le n' < N$ 
and $k \in \ZZ$, 
\[(\tfrac{N}{2}m+b(n'+kN),(n'+kN)\hspace*{-0.2cm} \mod N)
= (\tfrac{N}{2}(m+2bk)+bn',n'),\]
which implies
\[ 
\{(\tfrac{N}{2}m+bn,n \hspace*{-0.2cm} \mod N)\}_{m,
n\in\ZZ} =
\{(\tfrac{N}{2}m+bn,n)\}_{m\in\ZZ,
n=0,\ldots,N-1} \]
in the sense of sets. 

Using the definition of canonical generator matrices we can now define Gabor
systems for $\ltZ$.

\begin{definition}
Let $\Lambda$ be a lattice in $\ZZ \times \TT$ with canonical generator 
matrix $A$ given by
\begin{equation*}
A=\begin{bmatrix} \frac{N}{2} & b \\ 0 & \frac{1}{N} \end{bmatrix}
\end{equation*}
and let $g \in \ell_2(\ZZ)$. Then the associated Gabor 
system $\{g_{m\frac{N}{2}+nb,n\frac{1}{N}}\}_{m\in \ZZ,n=0,\dots,N-1}$ is given by
$$g_{m\frac{N}{2}+nb,n\frac{1}{N}}(l) = g(l-(m\tfrac{N}{2}+nb))e^{2 \pi i l n/N}, 
\qquad l\in \ZZ.$$
\end{definition}

Now we first give the definition of a Wilson basis associated with a lattice
with diagonal canonical generator matrix, i.e., with $b=0$ (in this
special case the definition coincides with the one given in~\cite{BFG97}).

\begin{definition}
\label{defiwil_rectangle}
Let $\Lambda$ be a lattice in $\ZZ \times \TT$ with canonical generator 
matrix $A$ given by
\begin{equation*}
A=\begin{bmatrix} \frac{N}{2} & 0 \\ 0 & \frac{1}{N} \end{bmatrix}
\end{equation*}
and let $g \in \ell_2(\ZZ)$.
Then the Wilson system 
$\Wil(g,\Lambda,\ell_2(\ZZ)) = \{\psi_{m,n}\}_{m\in\ZZ,
n=0,\ldots,\frac{N}{2}}$ 
is given by
\[ \hspace*{0.75cm} \psi_{m,n}^\Lambda = g_{mN,n\frac{1}{N}}, \hspace*{3.05cm} \mbox{if } 
m\in\ZZ,\;n=0,\tfrac{N}{2},\]
and for $m\in \ZZ$, $n=1,\ldots,\frac{N}{2}-1,$
\begin{eqnarray*}
\psi_{m,n}^\Lambda & = & \tfrac{1}{\sqrt{2}} 
(g_{m\frac{N}{2},n\frac{1}{N}}+ g_{m\frac{N}{2},-n\frac{1}{N}}), \quad
\mbox{if } m+n \mbox{ even},\\
\psi_{m,n}^\Lambda & = & \tfrac{i}{\sqrt{2}} 
(g_{m\frac{N}{2},n\frac{1}{N}}- g_{m\frac{N}{2},-n\frac{1}{N}}), \quad
\mbox{if } m+n \mbox{ odd}.
\end{eqnarray*}
\end{definition}

The Zak transform, which can be defined for any locally compact abelian
group (cf.~\cite{KK}), will be employed to prove equivalent conditions
for the Wilson system to form an orthonormal basis. In particular, we need
the Zak transform on $\TT$ with respect to the uniform lattice 
$K=\{e^{2 \pi i \frac{2k}{N}} : k=0,\ldots, \frac{N}{2}-1\}$ in $\TT$,
which is defined on the set of square--integrable functions on 
$\{e^{2 \pi i t}: t \in [0,\frac{2}{N})\} \times \{0,\ldots,\frac{N}{2}-1\}$ by
\[Zf(e^{2 \pi i t},y) = \sum_{k =0}^{\frac{N}{2}-1}f(e^{2 \pi i (t+\frac{2k}{N})})
e^{2 \pi i \frac{2k}{N}y}.\]

The proof of the following proposition is
inspired by the proof of \cite[Proposition 5.2]{DJJ91}.

\begin{proposition}
\label{eq_rec_Z}
Let $g \in \ell_2(\ZZ)$ be such that $\hat{g}$ is real--valued and consider the lattice $\Lambda$ 
with canonical generator matrix given by 
\begin{equation*}
\begin{bmatrix} \frac{N}{2} & 0 \\ 0 & \frac{1}{N} \end{bmatrix}.
\end{equation*}
Then the following conditions are equivalent.
\begin{enumerate}
\item[\rm (i)] $\{g_{m\frac{N}{2},n\frac{1}{N}}\}_{m\in\ZZ,n=0,\dots,N-1}$ 
is a tight frame for $\ell_2(\ZZ)$ with frame bound $2$.
\item[\rm (ii)] We have $|Z\hat{g}(e^{2 \pi i t},y)|^2 + |Z\hat{g}(e^{2 \pi i (t+\frac{1}{N})},y)|^2 
= N$ a.e..
\item[\rm (iii)] For all $j \in \{0,\ldots,N-1\}$, we have
$\sum_{l=0}^{N-1} \hat{g}(e^{2 \pi i (t+\frac{l}{N})}) 
\hat{g}(e^{2 \pi i (t+\frac{l+2j}{N})})
= N\delta_{j,0}$ a.e..
\item[\rm (iv)] $\Wil(g,\Lambda,\ell_2(\ZZ))$ is an orthonormal basis for $\ell_2(\ZZ)$.
\end{enumerate}
\end{proposition}

\begin{proof}
Throughout this proof we choose the normalized Haar measure on $\TT$, i.e.,
$m(E) = \int_0^1 1_{E}(e^{2 \pi i t}) dt$ for all measurable $E \subseteq \TT$
and the counting measure on its dual group $\widehat{\TT} = \ZZ$. This choice
ensures that the Plancherel formula for $\TT$ holds.

Since we will mainly work in the Fourier domain, we first need to compute the Fourier
transform of the elements of the Gabor system for our following calculations:
\begin{eqnarray*} 
\widehat{g_{m\frac{N}{2},n\frac{1}{N}}}(e^{2 \pi i t}) & = &
\sum_{l \in \ZZ} e^{2 \pi i \frac{n}{N}l} g(e^{2 \pi i (l-m\frac{N}{2})})
e^{-2 \pi i lt}\\
& = & e^{2 \pi i \frac{mn}{2}} e^{-2 \pi i tm\frac{N}{2}} 
\sum_{l \in \ZZ} g(e^{2 \pi i l}) e^{-2 \pi i l(t-\frac{n}{N})}\\
& = & (-1)^{mn} \hat{g}_{n\frac{1}{N},-m\frac{N}{2}}(e^{2 \pi i t}).
\end{eqnarray*}

(i) $\Leftrightarrow$ (ii):
Since the Fourier transform is a unitary operator, the Gabor system
$\{g_{m\frac{N}{2},\frac{1}{N}n}\}_{m\in\ZZ,n=0,\dots,N-1}$ 
is a tight frame for $\ell_2(\ZZ)$ with frame bound $2$ if and only if the Gabor system
$\{\hat{g}_{n\frac{1}{N},m\frac{N}{2}}\}_{m\in\ZZ,n=0,\dots,N-1}$ is a tight frame for 
$L^2(\TT)$ with frame bound $2$. Then we write this set as the disjoint union
$\{\hat{g}_{n\frac{2}{N},m\frac{N}{2}}\}_{m\in\ZZ,n=0,\dots,\frac{N}{2}-1} \cup
\{(T_{-\frac{1}{N}}\hat{g})_{n\frac{2}{N},m\frac{N}{2}}\}_{m\in\ZZ,n=0,\dots,\frac{N}{2}-1} 
=: G_1 \cup G_2$, and let $S_i$ denote the frame operator for $G_i$, $i=1,2$. 

First we compute the frame operator $S_1$. For all $f \in L^2(\TT)$, 
using the Poisson summation formula \cite[Theorem~4.42]{fol} applied to 
$H = \{e^{2 \pi i k\frac{2}{N}}: k = 0,\ldots,\frac{N}{2}-1\}$, we obtain
\allowdisplaybreaks
\begin{eqnarray*}
S_1f(e^{2 \pi i t})
& = & \sum_{m\in \ZZ}\sum_{n =0}^{\frac{N}{2}-1} \ip{f}{\hat{g}_{n\frac{2}{N},m\frac{N}{2}}}
\hat{g}_{n\frac{2}{N},m\frac{N}{2}}(e^{2 \pi i t})\\
& = & \sum_{m\in \ZZ}\sum_{n =0}^{\frac{N}{2}-1} \int_0^1 f(e^{2 \pi i s})
\overline{\hat{g}(e^{2 \pi i (s-n\frac{2}{N})})}e^{-2\pi i m\frac{N}{2}s} ds \:
\hat{g}(e^{2 \pi i (t-n\frac{2}{N})}) e^{2\pi i m\frac{N}{2}t}\\
& = & \sum_{n =0}^{\frac{N}{2}-1} \left[ \sum_{m\in \ZZ}
\widehat{(f\overline{T_{n\frac{2}{N}}\hat{g}})}(m\tfrac{N}{2})e^{2\pi i m\frac{N}{2}t}\right]
\hat{g}(e^{2 \pi i (t-n\frac{2}{N})})\\
& = & \sum_{n =0}^{\frac{N}{2}-1} \tfrac{2}{N}\sum_{k=0}^{\frac{N}{2}-1} 
(f\overline{T_{n\frac{2}{N}}\hat{g}})(e^{2 \pi i (t+k\frac{2}{N})})
\hat{g}(e^{2 \pi i (t-n\frac{2}{N})}).
\end{eqnarray*}
Applying now the Zak transform yields
\allowdisplaybreaks
\begin{eqnarray*}
Z(S_1f)(e^{2 \pi i t},y)
& = & \sum_{l =0}^{\frac{N}{2}-1}\tfrac{2}{N}\sum_{n,k =0}^{\frac{N}{2}-1}
f(e^{2 \pi i (t+\frac{2k+2l}{N})})\overline{\hat{g}(e^{2 \pi i (t+\frac{2k+2l-2n}{N})})}
\hat{g}(e^{2 \pi i (t+\frac{2l-2n}{N})}) e^{2 \pi i \frac{2l}{N}y}\\
& = & \tfrac{2}{N}\sum_{l,n,k =0}^{\frac{N}{2}-1}
f(e^{2 \pi i (t+\frac{2l+2n}{N})})\overline{\hat{g}(e^{2 \pi i (t+\frac{2l}{N})})}
\hat{g}(e^{2 \pi i (t+\frac{2l-2k}{N})}) e^{2 \pi i \frac{2(l-k+n)}{N}y}\\
& = & \tfrac{2}{N}\sum_{l =0}^{\frac{N}{2}-1}\overline{\hat{g}(e^{2 \pi i (t+\frac{2l}{N})})}
\left[ \sum_{k =0}^{\frac{N}{2}-1}\hat{g}(e^{2 \pi i (t+\frac{2l-2k}{N})})
e^{-2 \pi i \frac{2k}{N}y}\right]\\
&&\cdot\left[ \sum_{n =0}^{\frac{N}{2}-1}
f(e^{2 \pi i (t+\frac{2n+2l}{N})})e^{2 \pi i \frac{2n}{N}y}\right]
e^{2 \pi i \frac{2l}{N}y}\\
& = & \tfrac{2}{N}\overline{Z(\hat{g})(e^{2 \pi i t},y)}Z(\hat{g})(e^{2 \pi i t},y)
Z(f)(e^{2 \pi i t},y).
\end{eqnarray*}

To compute the Zak transform of $S_2$, we can use the
previous calculation, which yields
\begin{eqnarray*} 
Z(S_2f)(e^{2 \pi i t},y)& =& \tfrac{2}{N}Z(f)(e^{2 \pi i t},y)
|Z(T_{-\frac{1}{N}}\hat{g})(e^{2 \pi i t},y)|^2\\
& = & \tfrac{2}{N}Z(f)(e^{2 \pi i t},y)
|Z(\hat{g})(e^{2 \pi i (t+\frac{1}{N})},y)|^2,
\end{eqnarray*}
since
\[ Z(T_{-\frac{1}{N}}\hat{g})(e^{2 \pi i t},y) 
= \sum_{k =0}^{\frac{N}{2}-1}\hat{g}(e^{2 \pi i (t+\frac{2k}{N}+\frac{1}{N})})
e^{2 \pi i \frac{2k}{N}y}=Z(\hat{g})(e^{2 \pi i (t+\frac{1}{N})},y).\]
Hence (i) is equivalent to
\begin{eqnarray*} 
2Z(f)(e^{2 \pi i t},y)
& = & Z((S_1+S_2)f)(e^{2 \pi i t},y)\\
& = & \tfrac{2}{N}Z(f)(e^{2 \pi i t},y)\left[|Z(\hat{g})(e^{2 \pi i t},y)|^2 
+  |Z(\hat{g})(e^{2 \pi i (t+\frac{1}{N})},y)|^2\right] \mbox{ a.e.,}
\end{eqnarray*}
which holds if and only if (ii) is satisfied.

(ii) $\Leftrightarrow$ (iii):
The following properties of the Zak transform will be exploit several times.
The reconstruction formula
\[\sum_{y=0}^{\frac{N}{2}-1}Zf(e^{2 \pi i t},y)
=  \sum_{k =0}^{\frac{N}{2}-1}f(e^{2 \pi i (t+\frac{2k}{N})})
\sum_{y=0}^{\frac{N}{2}-1} e^{2 \pi i \frac{2k}{N}y}
= \tfrac{N}{2}f(e^{2 \pi i t})\]
holds a.e., since $\sum_{y=0}^{\frac{N}{2}-1} e^{2 \pi i \frac{2k}{N}y} \neq 0$ if and
only if $k=0$ by \cite[Lemma~23.29]{hew} and, if $k=0$, then 
$\sum_{y=0}^{\frac{N}{2}-1} e^{2 \pi i \frac{2k}{N}y} = \tfrac{N}{2}$. Moreover,
we will use that 
\[  Z\hat{g}(e^{2 \pi i (t+\frac{2l}{N})},y)
= \sum_{k =0}^{\frac{N}{2}-1}f(e^{2 \pi i (t+\frac{2k+2l}{N})})e^{2 \pi i \frac{2k}{N}y}
= e^{-2 \pi i \frac{2l}{N}y}Z\hat{g}(e^{2 \pi i t},y).\]

The idea is to write the equation in (iii) in terms of the Zak transform.
Using the fact that $\hat{g}$ is real--valued, we compute
\allowdisplaybreaks
\begin{eqnarray*}
\lefteqn{\sum_{l=0}^{N-1} \hat{g}(e^{2 \pi i (t+\frac{l}{N})}) 
\hat{g}(e^{2 \pi i (t+\frac{l+2j}{N})})}\\
& = & \frac{4}{N^2}\sum_{l=0}^{N-1} \sum_{x=0}^{\frac{N}{2}-1}Z\hat{g}(e^{2 \pi i (t+\frac{l}{N})},x)
\sum_{y=0}^{\frac{N}{2}-1}\overline{Z\hat{g}(e^{2 \pi i (t+\frac{l+2j}{N})},y)}\\
& = & \frac{4}{N^2}\sum_{k,x,y=0}^{\frac{N}{2}-1}\left[ Z\hat{g}(e^{2 \pi i (t+\frac{2k}{N})},x)
\overline{Z\hat{g}(e^{2 \pi i (t+\frac{2k+2j}{N})},y)}\right.\\
&&\left. +Z\hat{g}(e^{2 \pi i (t+\frac{2k+1}{N})},x)
\overline{Z\hat{g}(e^{2 \pi i (t+\frac{2k+2j+1}{N})},y)} \right]\\
& = & \frac{4}{N^2}\sum_{x,y=0}^{\frac{N}{2}-1}
\left[\sum_{k=0}^{\frac{N}{2}-1} e^{-2 \pi i \frac{2k}{N}(x-y)}\right]
e^{2 \pi i \frac{2j}{N}y}
\left[ Z\hat{g}(e^{2 \pi i t},x)
\overline{Z\hat{g}(e^{2 \pi i t},y)}\right.\\
&&\left. +Z\hat{g}(e^{2 \pi i (t+\frac{1}{N})},x)
\overline{Z\hat{g}(e^{2 \pi i (t+\frac{1}{N})},y)} \right]\\
& = & \frac{2}{N} \sum_{x=0}^{\frac{N}{2}-1}
\left[ |Z\hat{g}(e^{2 \pi i t},x)|^2
+|Z\hat{g}(e^{2 \pi i (t+\frac{1}{N})},x)|^2\right]e^{2 \pi i \frac{2j}{N}x}\\
& = & \frac{2}{N} \left[ |Z\hat{g}(e^{2 \pi i t},\cdot)|^2
+|Z\hat{g}(e^{2 \pi i (t+\frac{1}{N})},\cdot)|^2\right]^\vee(j),
\end{eqnarray*}
where the inverse Fourier transform is taken in $\ZZ_{\tfrac{N}{2}}$.
This shows that (iii) is equivalent to
\begin{equation}\label{ZakCL}\frac{2}{N} \left[ |Z\hat{g}(e^{2 \pi i t},\cdot)|^2
+|Z\hat{g}(e^{2 \pi i (t+\frac{1}{N})},\cdot)|^2\right]^\vee(j)  = N\delta_{j,0}.
\end{equation}
If (ii) holds, then
\[ \frac{2}{N} \left[ |Z\hat{g}(e^{2 \pi i t},\cdot)|^2
+|Z\hat{g}(e^{2 \pi i (t+\frac{1}{N})},\cdot)|^2\right]^\vee(j)
= \frac{2}{N} N \sum_{x=0}^{\frac{N}{2}-1}
 e^{2 \pi i \frac{2j}{N}x}
= N\delta_{j,0},\]
which is (\ref{ZakCL}). On the other hand, the inverse Fourier transform
is injective. This proves that (\ref{ZakCL}) holds if and only if 
(ii) is true and thus (ii) $\Leftrightarrow$ (iii).

(iii) $\Leftrightarrow$ (iv):
First we remark that $\Wil(g,\Lambda,\ell_2(\ZZ))$ is an orthonormal basis
if and only if the set
\[ \Psi := \{T_{nN}f_m : m=1,\ldots,N, \; n \in \ZZ\},\]
where
\begin{eqnarray*}
f_1(x) & = & g(x),\\
f_N(x) & = & g_{0,\frac{N}{2}}(x),\\
f_{2l+k}(x) & = & \frac{(-1)^{kl}}{\sqrt{2}} (g_{k\frac{N}{2},\frac{l}{N}} + (-1)^{k+l}
g_{k\frac{N}{2},-\frac{l}{N}}), \; l=1,\ldots,\tfrac{N}{2}-1,\: k=0,1
\end{eqnarray*}
is an orthonormal basis,
since these elements differ from the elements in $\Wil(g,\Lambda,\ell_2(\ZZ))$
only by factors of absolute value $1$. Next notice that to prove (iv) it is 
sufficient and necessary that
\begin{equation}\label{ONB1} \norm{T_{nN}f_m}_2=1, \quad m=1,\ldots,L, \; n \in \ZZ
\end{equation}
and
\begin{equation}\label{ONB2} \sum_{m=1}^L\sum_{n \in \ZZ} \ip{h_1}{T_{nN}f_m}\ip{T_{nN}f_m}{h_2}
=\ip{h_1}{h_2}, \quad  \mbox{for all }h_1,h_2 \in \ell_2(\ZZ).
\end{equation}
We start by dealing with (\ref{ONB1}). Using the Plancherel theorem, we compute
\[1 = \norm{T_{nN}f_1}_2^2 = \norm{\hat{g}}_2^2 = \int_0^1 \hat{g}(e^{2 \pi i t}) 
\hat{g}(e^{2 \pi i t}) dt,\]
\[1 = \norm{T_{nN}f_N}_2^2 = \|\hat{g}_{\frac{N}{2},0}\|^2 = \int_0^1 \hat{g}(e^{2 \pi i t}) 
\hat{g}(e^{2 \pi i t}) dt,\]
and, for $m=2,\ldots,N-1$, 
\begin{eqnarray*}
1 & = & \norm{T_{nN}f_m}_2^2\\
& = & \norm{\tfrac{1}{\sqrt{2}} (\hat{g}_{\frac{l}{N},-k\frac{N}{2}} + (-1)^{k+l}
\hat{g}_{-\frac{l}{N},-k\frac{N}{2}})}_2^2\\
& = & \frac12 \int_0^1 |e^{-2 \pi i \frac{kN}{2}t} |^2
|\hat{g}(e^{2 \pi i (t-\frac{l}{N})}) + (-1)^{k+l}\hat{g}(e^{2 \pi i (t+\frac{l}{N})})|^2 dt\\
& = & \frac12 \int_0^1 \left[ |\hat{g}(e^{2 \pi i (t-\frac{l}{N})})|^2 + |\hat{g}(e^{2 \pi i
(t+\frac{l}{N})})|^2 + (-1)^{k+l}\hat{g}(e^{2 \pi i (t-\frac{l}{N})})\overline{\hat{g}(e^{2
\pi i (t+\frac{l}{N})})}\right.\\
&& \left.+ (-1)^{k+l}\overline{\hat{g}(e^{2 \pi i (t-\frac{l}{N})})}
\hat{g}(e^{2 \pi i (t+\frac{l}{N})})\right] dt.
\end{eqnarray*}
Since $\hat{g}$ is real--valued, we can continue the last computation and obtain that
\[ 1 = \norm{T_{nN}f_m}_2^2 = \norm{\hat{g}}_2^2 + (-1)^{k+l} \int_0^1
\hat{g}(e^{2 \pi i t}) \hat{g}(e^{2 \pi i (t+\frac{2l}{N})}) dt.\]
Combining the above computations we have proven that (\ref{ONB1}) holds
if and only if 
\begin{equation}\label{ONB1_eq} \int_0^1
\hat{g}(e^{2 \pi i t}) \hat{g}(e^{2 \pi i (t+\frac{2j}{N})}) dt = \delta_{j,0} \quad
\mbox{for all }j \in \{0,\ldots,N-1\}.
\end{equation}

Now we turn to the study of condition (\ref{ONB2}).
Using the Plancherel formula and the Poisson summation formula 
\cite[Theorem~4.42]{fol} applied to $H = \{e^{2 \pi i \frac{k}{N}}: k 
= 0,\ldots,N-1\}$, we obtain
\allowdisplaybreaks
\begin{eqnarray*}
\lefteqn{\sum_{m=1}^N\sum_{n \in \ZZ} \ip{h_1}{T_{nN}f_m}\ip{T_{nN}f_m}{h_2}}\\
& = & \sum_{m=1}^N\sum_{n \in \ZZ} \ip{\widehat{h_1}}{\widehat{T_{nN}f_m}}
\ip{\widehat{T_{nN}f_m}}{\widehat{h_2}}\\
& = & \sum_{m=1}^N\sum_{n \in \ZZ} \int_0^1 
(\widehat{h_1}\overline{\widehat{f_m}})(e^{2 \pi i t}) \int_0^1 
(\widehat{f_m}\overline{\widehat{h_2}})(e^{2 \pi i s}) e^{-2 \pi i sNn}ds \:
e^{2 \pi i tNn} dt\\
& = & \sum_{m=1}^N \int_0^1 (\widehat{h_1}\overline{\widehat{f_m}})(e^{2 \pi i t}) \left[
\sum_{n \in \ZZ}(\widehat{f_m}\overline{\widehat{h_2}})\widehat{\;\:}(Nn)e^{2 \pi i tNn}
\right] dt\\
& = & \sum_{m=1}^N \int_0^1 (\widehat{h_1}\overline{\widehat{f_m}})(e^{2 \pi i t})
\tfrac{1}{N} 
\sum_{r =0}^{N-1}(\widehat{f_m}\overline{\widehat{h_2}})(e^{2 \pi i (t+\frac{r}{N})}) \: dt,
\end{eqnarray*}
which equals $\ip{h_1}{h_2}$ if and only if
\begin{equation}\label{ONB2_2}
\sum_{m=1}^N \overline{\widehat{f_m}(e^{2 \pi i t})} 
\widehat{f_m}(e^{2 \pi i (t+\frac{r}{N})}) = N\delta_{r,0}
\quad \mbox{for all }r \in \{0,\ldots,N-1\}.
\end{equation}
Setting $\mathbb{L} := \{-\tfrac{N}{2}+1,\ldots,-1,1,\ldots,\tfrac{N}{2}-1\}$, 
we compute
\allowdisplaybreaks
\begin{eqnarray*}
\lefteqn{\sum_{m=1}^N \overline{\widehat{f_m}(e^{2 \pi i t})}
\widehat{f_m}(e^{2 \pi i (t+\frac{r}{N})})}\\
& = & \hat{g}(e^{2 \pi i t})\hat{g}(e^{2 \pi i (t+\frac{r}{N})}) 
+ \hat{g}_{\frac{N}{2},0}(e^{2 \pi i t})\hat{g}_{\frac{N}{2},0}(e^{2 \pi i (t+\frac{r}{N})})
+ \frac12 \sum_{l=1}^{\frac{N}{2}-1}\sum_{k=0}^1 
[\overline{\hat{g}_{\frac{l}{N},-k\frac{N}{2}}(e^{2 \pi i t})}\\
&& + (-1)^{k+l} \overline{\hat{g}_{-\frac{l}{N},-k\frac{N}{2}}(e^{2 \pi i t})}]
[\hat{g}_{\frac{l}{N},-k\frac{N}{2}}(e^{2 \pi i (t+\frac{r}{N})}) + 
(-1)^{k+l} \hat{g}_{-\frac{l}{N},-k\frac{N}{2}}(e^{2 \pi i (t+\frac{r}{N})})]\\
& = & \hat{g}(e^{2 \pi i t})\hat{g}(e^{2 \pi i (t+\frac{r}{N})}) 
+\hat{g}(e^{2 \pi i (t-\frac{N}{2})})\hat{g}(e^{2 \pi i (t+\frac{r}{N}-\frac{N}{2})}) 
+ \frac12 \sum_{l=1}^{\frac{N}{2}-1}\sum_{k=0}^1 e^{-2 \pi i k\frac{N}{2}\frac{r}{N}}
[\hat{g}(e^{2 \pi i (t-\frac{l}{N})})\\
&&+(-1)^{k+l}\hat{g}(e^{2 \pi i (t+\frac{l}{N})})]
[\hat{g}(e^{2 \pi i (t-\frac{l}{N}+\frac{r}{N})})+
(-1)^{k+l}\hat{g}(e^{2 \pi i (t+\frac{l}{N}+\frac{r}{N})})]\\
& = & \hat{g}(e^{2 \pi i t})\hat{g}(e^{2 \pi i (t+\frac{r}{N})}) 
+\hat{g}(e^{2 \pi i (t-\frac{N}{2})})\hat{g}(e^{2 \pi i (t+\frac{r}{N}-\frac{N}{2})}) \\
&& + \frac12 \sum_{l=1}^{\frac{N}{2}-1}\sum_{k=0}^1  (-1)^{kr}
\left[\hat{g}(e^{2 \pi i (t-\frac{l}{N})})\hat{g}(e^{2 \pi i (t-\frac{l}{N}+\frac{r}{N})})
+\hat{g}(e^{2 \pi i (t+\frac{l}{N})})\hat{g}(e^{2 \pi i (t+\frac{l}{N}+\frac{r}{N})})
\right.\\
&& \left. + (-1)^{k+l}
(e^{2 \pi i (t-\frac{l}{N})})\hat{g}(e^{2 \pi i (t+\frac{l}{N}+\frac{r}{N})})
+\hat{g}(e^{2 \pi i (t+\frac{l}{N})})\hat{g}(e^{2 \pi i
(t-\frac{l}{N}+\frac{r}{N})}))\right]\\
& = & \hat{g}(e^{2 \pi i t})\hat{g}(e^{2 \pi i (t+\frac{r}{N})})  
+ \sum_{l \in \mathbb{L}}\hat{g}(e^{2 \pi i (t+\frac{l}{N})})\hat{g}(e^{2 \pi i
(t+\frac{l}{N}+\frac{r}{N})})\left[\frac{1+(-1)^r}{2}\right]\\
&&+\sum_{l \in \mathbb{L}}(-1)^l \hat{g}(e^{2 \pi i (t+\frac{l}{N})})\hat{g}(e^{2 \pi i
(t-\frac{l}{N}+\frac{r}{N})})\left[\frac{1+(-1)^{r+1}}{2}\right]
+\hat{g}(e^{2 \pi i (t-\frac{N}{2})})\hat{g}(e^{2 \pi i (t+\frac{r}{N}-\frac{N}{2})}).
\end{eqnarray*}
If $r$ is even, i.e., $r=2j$, we obtain
\[\sum_{m=1}^N \overline{\widehat{f_m}(e^{2 \pi i t})}
\widehat{f_m}(e^{2 \pi i (t+\frac{r}{N})})
= \sum_{l=0}^{N-1} \hat{g}(e^{2 \pi i (t+\frac{l}{N})})
\hat{g}(e^{2 \pi i (t+\frac{l}{N}+\frac{2j}{N})}),\]
and if $r$ is odd, i.e., $r=2j+1$, we obtain
\begin{eqnarray*} 
\lefteqn{\sum_{m=1}^N \overline{\widehat{f_m}(e^{2 \pi i t})}
\widehat{f_m}(e^{2 \pi i (t+\frac{r}{N})})}\\
& = & \hat{g}(e^{2 \pi i t})\hat{g}(e^{2 \pi i (t+\frac{r}{N})}) 
+\hat{g}(e^{2 \pi i (t-\frac{N}{2})})\hat{g}(e^{2 \pi i (t+\frac{r}{N}-\frac{N}{2})})\\
&&
+ \sum_{l \in \mathbb{L}}(-1)^l \hat{g}(e^{2 \pi i (t+\frac{l}{N})})\hat{g}(e^{2 \pi i
(t-\frac{l}{N}+\frac{r}{N})})\\
& =&0.
\end{eqnarray*}
This shows that (\ref{ONB2}) holds if and only if
\[\sum_{l=0}^{N-1} \hat{g}(e^{2 \pi i (t+\frac{l}{N})}) 
\hat{g}(e^{2 \pi i (t+\frac{l+2j}{N})})
= N\delta_{j,0}.\]
Moreover, this equation implies equation (\ref{ONB1_eq}), since
\[ \int_0^1 \hat{g}(e^{2 \pi i t}) \hat{g}(e^{2 \pi i (t+\frac{2j}{N})}) dt
= \int_0^{\frac{1}{N}} \sum_{l=0}^{N-1}
\hat{g}(e^{2 \pi i (t+\frac{l}{N})}) \hat{g}(e^{2 \pi i (t+\frac{l+2j}{N})}) dt
= \int_0^{\frac{1}{N}} N\delta_{j,0} dt = \delta_{j,0}.\]
This shows (iii) $\Leftrightarrow$ (iv), and hence the theorem is proven.
\end{proof}

Now we will study the case of a general time--frequency lattice. 

\begin{proposition}
\label{meta_Z}
Let $\Lambda$ be a lattice in $\ZZ \times \TT$ with $\vol(\Lambda) = \frac12$ with
canonical generator matrix $A$ given by
\begin{equation*}
A=\begin{bmatrix}
\frac{N}{2} & b \\
0 & \frac{1}{N}
\end{bmatrix}.
\end{equation*}
Let $g \in \ell_2(\ZZ)$, let $m_0,n_0 \in \ZZ$ be chosen
such that $\frac{N}{2}m_0+bn_0=\gcd(\frac{N}{2},b)=:c$, and  let $U$ be defined 
on $\ell_2(\ZZ)$ by
\[ Uf(k) = f(k)e^{\pi i \frac{n_0}{cN}k^2}.\]
Then 
\[ g_{m\frac{N}{2}+nb,n\frac{1}{N}}(l) =
C(m,n) U(U^{-1}g)_{m\frac{N}{2}+nb,-m\frac{n_0}{2c}-n\frac{bn_0}{cN}+n\frac{1}{N}}(l),
\quad l \in \ZZ,\]
where $C(m,n) = e^{\pi i \frac{n_0}{cN}(m\frac{N}{2}+nb)^2}$.
\end{proposition}

\begin{proof}
Let $\sigma \in {\rm Hom}(\ZZ \times \TT)$ be defined by
\[ \sigma = \begin{bmatrix} I_\ZZ & 0 \\ -\frac{n_0}{cN} & I_\TT \end{bmatrix}.\]
It is easy to check that $\sigma$ is symplectic on $\ZZ \times \TT$. In order to apply
Theorem \ref{Kai_Theorem}, we need to compute a second degree character
of $\ZZ \times \TT$ associated to $\zeta = \sigma^*\kappa_0\sigma-\kappa_0$,
where $\kappa_0$ is defined by 
\[ \kappa_0 = \begin{bmatrix} 1 & 1 \\ I_{\ZZ} & 0 \end{bmatrix}
\in {\rm Hom}(\ZZ \times \TT,\TT \times \ZZ).\]
First we note that $\sigma^*(z,n) = (ze^{2 \pi i (-\frac{n_0}{cN})n},n)$, since
\[ \ip{(m,t)}{\sigma^*(z,n)} 
= \ip{\sigma(m,t)}{(z,n)} 
= z^m (e^{2 \pi i (-\frac{n_0}{cN})m} t)^n
= \ip{(m,t)}{(ze^{2 \pi i (-\frac{n_0}{cN})n},n)}.\]
Thus
\[\zeta(n,z) 
=  (\sigma^*\kappa_0\sigma-\kappa_0)(n,z)
= \sigma^*\kappa_0(n,e^{2 \pi i (-\frac{n_0}{cN})n}z)-(1,n)
= (e^{2 \pi i (-\frac{n_0}{cN})n},0).\]
The map 
\[ \psi : \ZZ \times \TT \to \TT,\quad \psi(m,t) = e^{-\pi i \frac{n_0}{cN} m^2}\]
is a second degree character associated to $\zeta$ as the following
calculation shows:
\begin{eqnarray*} 
\psi(m,t)\psi(n,z)\ip{(m,t)}{\zeta(n,z)}
& = & e^{-\pi i \frac{n_0}{cN} m^2} e^{-\pi i \frac{n_0}{cN} n^2}
e^{2 \pi i (-\frac{n_0}{cN})mn}\\
& = & e^{-\pi i \frac{n_0}{cN} (m+n)^2}\\
& = & \psi((m,t)+(n,z)).
\end{eqnarray*}
Next notice that
\begin{equation}\label{sigma_connection}
\sigma(m\tfrac{N}{2}+nb,n\tfrac{1}{N})
= (m\tfrac{N}{2}+nb,-m\tfrac{n_0}{2c}-n\tfrac{bn_0}{cN}+n\tfrac{1}{N}).
\end{equation}
Now we can apply Theorem \ref{Kai_Theorem}, 
which proves the claim.
\end{proof}

Next we define a Wilson basis associated with a lattice
with arbitrary canonical generator matrix. For this, 
the following mapping will turn out to be very useful.

\begin{lemma}
\label{def_varphi}
Let $\frac{N}{2}, b \in \ZZ$ with $0 \le b < \frac{N}{2}$, and
let $m_0,n_0 \in \ZZ$ be chosen
such that $\frac{N}{2}m_0+bn_0=\gcd(\frac{N}{2},b)=:c$.
Further let $d:= {\rm lcm}(\frac{N}{2},b)$. Then
the mapping $\varphi : \ZZ^2 \to \ZZ^2$ defined by
\[ \phi(m,n) = \left\{\begin{array}{ccl}
(m,n) & : & b=0,\\
(mm_0-\frac{2d}{N}n,mn_0+\frac{d}{b}n) & : & b \neq 0 
\end{array}\right.\]
is bijective and, for all $m \in \ZZ$, we have
\[
\{(m,n \hspace*{-0.2cm} \mod 2c) : (m,n) \in \varphi^{-1}(\ZZ 
\times \{0,\ldots,N-1\})\} = \{m\} \times \{0,\ldots,2c-1\}
\]
with
\[|\{n : (m,n) \in \varphi^{-1}(\ZZ 
\times \{0,\ldots,N-1\})\}|=2c.\]
\end{lemma}

\begin{proof}
We only need to study the case $b \neq 0$. For this, 
let $(m,n),(m',n')\in\ZZ^2$ be such that $\phi(m,n) = \phi(m',n')$.
Then
\[ \tfrac{N}{2}(mm_0-\tfrac{2d}{N}n) + b(mn_0+\tfrac{d}{b}n) 
= \tfrac{N}{2}(m'm_0-\tfrac{2d}{N}n') + b(m'n_0+\tfrac{d}{b}n'),\]
which holds if and only if
\[ m(\tfrac{N}{2}m_0+bn_0) = m'(\tfrac{N}{2}m_0+bn_0),\]
and hence $m=m'$. This implies
\[ (-\tfrac{2d}{N}n,\tfrac{d}{b}n) = (-\tfrac{2d}{N}n',\tfrac{d}{b}n'),\]
which yields $n=n'$. This proves that $\varphi$ is injective.

To show that $\varphi$ is surjective, let $(k,l) \in \ZZ^2$ and
consider $M:=\frac{N}{2}k+bl$. It is well--known that
there exists some $m \in \ZZ$ with $M=mc$. Furthermore, we
have
\[ \{(p,q)\in \ZZ^2 : \tfrac{N}{2}p+bq=mc\}
= \{(mm_0-\tfrac{2d}{N}n,mn_0+\tfrac{d}{b}n):n \in \ZZ\},\]
since $\tfrac{N}{2}p+bq=\frac{N}{2}p'+bq'$ if and only if
$\frac{N}{2}(p-p') = b(q'-q)$.
This yields the existence of some $n \in \ZZ$ with 
\[\varphi(m,n)=(mm_0-\tfrac{2d}{N}n,mn_0+\tfrac{d}{b}n) = (k,l).\]

Secondly, we will prove the second part of the lemma. First observe that
$m,n \in \ZZ$ satisfy 
\begin{equation} \label{eq_proof_lem} \varphi(m,n)
\in \ZZ \times \{0,\ldots,N-1\}
\end{equation} if and only if they satisfy
\[ -\tfrac{b}{d}n_0m \le n \le -\tfrac{b}{d}n_0m + \tfrac{b}{d}(N-1)
= 2c-\tfrac{b}{d}n_0m - \tfrac{b}{d}.\]
Hence, for each fixed $m \in \ZZ$, the set of $n \in \ZZ$ such that
(\ref{eq_proof_lem}) is satisfied
equals 
\[ S_m := \{\lceil -\tfrac{b}{d}n_0m \rceil,\ldots,\lfloor 
2c-\tfrac{b}{d}n_0m - \tfrac{b}{d} \rfloor\}.\]

To finish the proof we claim that 
\begin{equation}\label{eq_proof_lem_2}
|S_m|=2c \quad \mbox{ for all } m \in \ZZ.\end{equation}
For this, fix $m \in \ZZ$ and let $k \in \ZZ$ and $l \in \{0,\ldots,\frac{d}{b}-1\}$ 
be such that $-n_0m = k \frac{d}{b} + l$. Then we obtain
$\lceil -\tfrac{b}{d}n_0m \rceil = k$ if $l = 0$ and otherwise 
$\lceil -\tfrac{b}{d}n_0m \rceil =k+1$. Moreover, we have
\[ \lfloor  2c-\tfrac{b}{d}n_0m - \tfrac{b}{d} \rfloor
= \left\lfloor  2c + k + \tfrac{l-1}{\frac{d}{b}} \right\rfloor,\]
which equals $2c+k-1$ if $l=0$ and otherwise $2c+k$. Thus
the second part of the lemma is proven.
\end{proof}

Note that the following definition reduces to Definition \ref{defiwil_rectangle} 
in the case of a diagonal canonical generator matrix. 

\begin{definition}
Let $\Lambda$ be a lattice in $\ZZ \times \TT$ with canonical generator 
matrix $A$ given by
\begin{equation*}
A=\begin{bmatrix} \frac{N}{2} & b \\ 0 & \frac{1}{N} \end{bmatrix}.
\end{equation*}
Let $g \in \ell_2(\ZZ)$, and let
$\varphi$ be defined as in  Lemma \ref{def_varphi}.
Then the Wilson system 
$\Wil(g,\Lambda,\ell_2(\ZZ)) = \{\psi_{m,n}\}_{m \in \ZZ, n = 0,\ldots,\frac{N}{2}}$ 
is given by
\[ \hspace*{0.85cm} \psi_{m,n}^\Lambda = 
g_{\varphi_1(2m,n)\frac{N}{2}+\varphi_2(2m,n)b,\varphi_2(2m,n)\frac{1}{N}}, 
\hspace*{0.9cm} \mbox{if } m\in\ZZ,\;n=0,\tfrac{N}{2},\]
and for $m\in \ZZ$, $n=1,\ldots,\frac{N}{2}-1,$
\begin{eqnarray*}
\psi_{m,n}^\Lambda & = & \tfrac{1}{\sqrt{2}} 
(g_{\varphi_1(m,n)\frac{N}{2}+\varphi_2(m,n)b,\varphi_2(m,n)\frac{1}{N}}\\
&&
+ g_{\varphi_1(m,-n)\frac{N}{2}+\varphi_2(m,-n)b,\varphi_2(m,-n)\frac{1}{N}}), \quad
\mbox{if } m+n \mbox{ even},\\
\psi_{m,n}^\Lambda & = & \tfrac{i}{\sqrt{2}} 
(g_{\varphi_1(m,n)\frac{N}{2}+\varphi_2(m,n)b,\varphi_2(m,n)\frac{1}{N}}\\
&&
- g_{\varphi_1(m,-n)\frac{N}{2}+\varphi_2(m,-n)b,\varphi_2(m,-n)\frac{1}{N}}), \quad
\mbox{if } m+n \mbox{ odd}.
\end{eqnarray*}
\end{definition}

The following theorem gives an equivalent condition for a Wilson
system with respect to an arbitrary time--frequency lattice to
form an orthonormal basis in terms of a frame condition for the
associated Gabor system.

\begin{theorem}
\label{theo_general_Z}
Let $\Lambda$ be a lattice in $\ZZ \times \TT$ with canonical generator 
matrix $A$ given by
\begin{equation*}
A=\begin{bmatrix} \frac{N}{2} & b \\ 0 & \frac{1}{N} \end{bmatrix}.
\end{equation*}
Let $g \in \ell_2(\ZZ)$ be such that  $\widehat{U^{-1}g}$ is real--valued, let $M:=2c$, 
and let $U$ and $\varphi$ be defined as in
Proposition \ref{meta_Z} and Lemma \ref{def_varphi}, respectively.
Then the following conditions are equivalent.
\begin{enumerate}
\item[\rm (i)] $\{g_{m\frac{N}{2}+nb,n\frac{1}{N}}\}_{m\in \ZZ,n=0,\dots,N-1}$ is a tight 
frame for $\ell_2(\ZZ)$ with frame bound $2$.
\item[\rm (ii)] $\{(U^{-1}g)_{m\frac{M}{2},n\frac{1}{M}}\}_{m\in \ZZ,n=0,\dots,M-1}$ 
is a tight frame for $\ell_2(\ZZ)$ with frame bound $2$.
\item[\rm (iii)] $\Wil(U^{-1}g,\frac{M}{2}\ZZ \times \frac{1}{M}\{0,\ldots,M-1\},
\ell_2(\ZZ))$ is an orthonormal basis for $\ell_2(\ZZ)$.
\item[\rm (iv)] $\Wil(g,\Lambda,\ell_2(\ZZ))$ is an orthonormal basis for $\ell_2(\ZZ)$.
\end{enumerate}
\end{theorem}

\begin{proof}
Let $\sigma$ be defined as in the proof of Proposition \ref{meta_Z} and
let $\varphi = (\varphi_1,\varphi_2)$.
Then we compute
\begin{eqnarray*}
\lefteqn{\sigma(\varphi_1(m,n)\tfrac{N}{2}+\varphi_2(m,n)b,\varphi_2(m,n)\tfrac{1}{N})}\\
& = & (\varphi_1(m,n)\tfrac{N}{2}+\varphi_2(m,n)b,
-\varphi_1(m,n)\tfrac{n_0}{2c}-\varphi_2(m,n)\tfrac{bn_0}{cN}+ \varphi_2(m,n)\tfrac{1}{N})\\
& = & m(\tfrac{N}{2}m_0+bn_0)+n(-\tfrac{2d}{N}\tfrac{N}{2}+\tfrac{d}{b}b),
m(-\tfrac{n_0}{2c}m_0-\tfrac{bn_0}{cN}n_0+\tfrac{1}{N}n_0)\\
&& +n(\tfrac{2d}{N}\tfrac{n_0}{2c}
-\tfrac{d}{b}\tfrac{bn_0}{cN}+\tfrac{d}{b}\tfrac{1}{N}))\\
& = & (mc,mn_0(\tfrac{1}{N}-\tfrac{1}{cN}(\tfrac{N}{2}m_0+bn_0))
+n\tfrac{d}{bN})\\
& = & (m\tfrac{M}{2},n\tfrac{1}{M}),
\end{eqnarray*}
where in the last step we used $cd=\frac{N}{2} b$.
Using Lemma \ref{def_varphi}, the equivalence of 
(i) and (ii) now follows immediately from Proposition \ref{meta_Z} and
Equation (\ref{sigma_connection}), since $U$ is unitary and $|C(m,n)|=1$.
Proposition \ref{eq_rec_Z} proves (ii) $\Leftrightarrow$ (iii). 
Therefore it remains to prove the equivalence of (iii) and (iv). 
For this, we will use the following 
implication of Proposition \ref{meta_Z}:
\begin{eqnarray*} 
U(U^{-1}g)_{m\frac{M}{2},n\frac{1}{M}}
&=& U(U^{-1}g)_{\sigma(\varphi_1(m,n)\frac{N}{2}+\varphi_2(m,n)b,\varphi_2(m,n)\frac{1}{N})}\\
&=& C(\varphi(m,n))^{-1}g_{\varphi_1(m,n)\frac{N}{2}+\varphi_2(m,n)b,\varphi_2(m,n)\frac{1}{N}}.
\end{eqnarray*}
Further notice that $C(\varphi(m,n))^{-1}$ does not depend on the sign of $n$, since
\[   C(\varphi(m,n))^{-1} = e^{-\pi i \frac{n_0}{cN}(\varphi_1(m,n)\frac{N}{2}+\varphi_2(m,n)b)^2}
= e^{-\pi i \frac{n_0}{cN}m^2(m_0\frac{N}{2}+n_0b)^2}
= e^{-\pi i \frac{n_0}{N}m^2c^2}.\]
Using now the definition of a Wilson basis, the fact that $U$ is a unitary
operator, and the fact that $|C(\varphi(m,n))|=1$ yields the result. 
\end{proof}


\section{Wilson bases for general lattices -- the finite case} 
\label{s:fin}

The space $\Cst^L$ has several advantages over $\ltZ$ when
constructing numerical methods for practical
time-frequency analysis, which often allow a further acceleration of 
numerical algorithms, e.g., see~\cite{AT90,Str97a}.

Before defining Gabor systems and Wilson systems for $\Cst^L$ for general 
time-frequency lattices, 
we first prove that each such lattice does not only possess a uniquely 
determined generator matrix in Hermite normal form (which was already proved 
in \cite{Her1850}), but moreover in our situation this matrix attains a
special form. 

\begin{lemma}
\label{latred2}
Let $\Lambda$ be a lattice in $\ZZ_L \times \ZZ_L$ with generator matrix $A$ given by
\begin{equation*}
A=\begin{bmatrix} a & b \\ c & d \end{bmatrix},
\qquad \text{with $a,b,c,d \in \Nst$, and $\det(A)=\frac{L}{2}$},
\end{equation*}
and denote $p=\gcd(c,d)$ if $c\neq 0$ and $p=d$ if $c=0$. 
Then $\Lambda$ possesses a uniquely determined generator matrix of 
the form
\begin{equation}
A' =
\begin{bmatrix}
 \frac{L}{2p} & b' \\
0 & p
\end{bmatrix},
\label{Apmat}
\end{equation}
where $p=\gcd(c,d)$ and $0 \le b' < \frac{L}{2p}$.
\end{lemma}

\begin{proof}
Let $p = \gcd(c,d)$ with $c=qp, d=rp$ and note that $p | \Lt$ since
$\det(A) = ad-bc=(ar-bq)p= \Lt$. Since $\frac{d}{p}=r$ and $\frac{-c}{p}=-q$,
we have
\begin{equation}
A \begin{bmatrix} r \\ -q \end{bmatrix} =
\begin{bmatrix} \frac{L}{2p} \\ 0 \end{bmatrix}.
\label{Amat1}
\end{equation}
Furthermore, we claim that there exists a $z \in \ZZ$ with $0 \le z < \frac{L}{2p}$
such that the point
$\begin{bmatrix} z \\ p \end{bmatrix}$
belongs to $\Lambda$. This can be seen as follows: The condition
$\begin{bmatrix} z \\ p \end{bmatrix} \in \Lambda$ is equivalent to the existence
of $m,n \in \Zst$ such that
\begin{equation}
\begin{bmatrix} a & b \\ c & d \end{bmatrix}
 \begin{bmatrix} m \\ n \end{bmatrix} =
\begin{bmatrix} z \\ p \end{bmatrix}.
\label{Amat2}
\end{equation}
Consider the equation $cm+dn=p$ and substitute $c=qp, d=rp$, then
$qpm + rpn = p$, hence $qm+rn=1$. Since $r$ and $q$ are relative prime,
there exist $m,n \in \Zst$ such that $qm+rn=1$
(see e.g. \cite[Theorem~4.4]{Hua82}). Thus~\eqref{Amat2} holds for
$z \in \Zst$, but we still have to show that it holds under the condition
$0\le z < \frac{L}{2p}$. We can write $z=b'+k\frac{L}{2p}$ with
$0\le b' < \frac{L}{2p}$ and $k\in \Zst$.  Hence
$$
\begin{bmatrix} z \\ p \end{bmatrix}=
\begin{bmatrix} b' \\ p \end{bmatrix}+ k
\begin{bmatrix} \frac{L}{2p} \\ 0 \end{bmatrix}.
$$
Since $ \begin{bmatrix} \frac{L}{2p} \\ 0 \end{bmatrix} =  A
\begin{bmatrix} r \\ -q \end{bmatrix}$ by~\eqref{Amat1}, it follows that
$\begin{bmatrix} b' \\ p \end{bmatrix} \in \Lambda$.
Consequently the matrix
\begin{equation*}
A' =
\begin{bmatrix}
\frac{L}{2p} & b' \\
0 & p
\end{bmatrix}
\end{equation*}
(which satisfies $\det(A') = \Lt$) generates $\Lambda$.

The fact that this matrix is uniquely determined is an immediate consequence from
the condition $0 \le b' < \frac{L}{2p}$.
\end{proof}

\begin{definition}
Let $\Lambda$ be a lattice in $\ZZ_L \times \ZZ_L$. Then the uniquely determined
matrix $A'$ of Lemma \ref{latred2} is called the
{\em canonical generator matrix} for $\Lambda$. 
\end{definition}

Using the notion of a canonical generator matrix, we first give the definition
of a Gabor system.

\begin{definition}
Let $\Lambda$ be a lattice in $\ZZ_L \times \ZZ_L$ with canonical generator 
matrix $A$ given by
\begin{equation*}
A=\begin{bmatrix} \frac{L}{2p} & b \\ 0 & p \end{bmatrix}.
\end{equation*}
Set $M=2p, N=\frac{L}{p}$ and let $g$ be some $L$-periodic function 
on $\Zst$. Then the associated Gabor 
system is given by $\{g_{ma+nb,nd}\}_{m=0,\dots,M-1,n=0,\dots,N-1}$, where
$$g_{ma+nb,nd}(l) = g(l-(ma+nb))e^{2 \pi i l nd/L}, 
\qquad l=0,\dots,L-1.$$
\end{definition}

Next we define a Wilson basis associated with a lattice
with diagonal canonical generator matrix in the following way:

\begin{definition}
\label{defiwil_rectangleCL}
Let $\Lambda$ be a lattice in $\ZZ_L \times \ZZ_L$ with canonical generator 
matrix $A$ given by
\begin{equation*}
A=\begin{bmatrix} \frac{L}{2p} & 0 \\ 0 & p \end{bmatrix}
\end{equation*}
and let $g$ be some $L$-periodic function 
on $\Zst$.
Then the Wilson system 
$\Wil(g,\Lambda,\Cst^L) = \{\psi_{m,n}\}_{(m,n) \in I}$, where $I = \{0,\ldots,p-1\}
\times \{0,\frac{L}{2p}\} \cup \{0,\ldots,2p-1\} \times \{1,\ldots,\frac{L}{2p}-1\}$,
is given by
\[ \hspace*{2.3cm}\psi_{m,n}^\Lambda = g_{m\frac{L}{p},np}, \hspace*{2.9cm}\mbox{if }
m=0,\ldots,p-1,\; n=0,\tfrac{L}{2p},\]
and for $m=0,\ldots,2p-1$, $n=1,\ldots,\frac{L}{2p}-1,$
\begin{eqnarray*}
\psi_{m,n}^\Lambda & = & \tfrac{1}{\sqrt{2}} 
(g_{m\frac{L}{2p},np}+ g_{m\frac{L}{2p},-np}), \quad
\mbox{if } m+n \mbox{ even},\\
\psi_{m,n}^\Lambda & = & \tfrac{i}{\sqrt{2}} 
(g_{m\frac{L}{2p},np}- g_{m\frac{L}{2p},-np}), \quad
\mbox{if } m+n \mbox{ odd}.
\end{eqnarray*}
\end{definition}

Also in the finite case we will employ the Zak transform. This time
we will use the Zak transform on the group $\ZZ_L$ with respect to 
the uniform lattice $K=\{2pk: k=0,\ldots, \frac{L}{2p}-1\}$ in $\ZZ_L$, which 
is defined on the set of square--integrable functions on the set 
$\{0,\ldots,2p-1\} \times \{0,\ldots,\frac{L}{2p}-1\}$ by
\[Zf(x,y) = \sum_{k =0}^{\frac{L}{2p}-1}f(x+2pk)
e^{2 \pi i \frac{2pk}{L}y},\]
where we associate $\ZZ_L$ with $\{0,\ldots,L-1\}$.

The following proposition is the analog to Proposition~\ref{eq_rec_Z} for
the space $\CC^L$.
\begin{proposition}
\label{eq_rec_CL}
Let $g$ be some $L$--periodic function on $\ZZ$ such that $\hat{g}$ is real--valued 
and consider the lattice $\Lambda$ 
with canonical generator matrix given by 
\begin{equation*}
\begin{bmatrix} \frac{L}{2p} & 0 \\ 0 & p \end{bmatrix}.
\end{equation*}
Then the following conditions are equivalent.
\begin{enumerate}
\item[\rm (i)] $\{g_{m\frac{L}{2p},np}\}_{m=0,\ldots,2p-1,n=0,\dots,\frac{L}{p}-1}$ is a tight frame for $\CC^L$
with frame bound $2$.
\item[\rm (ii)] We have $|Z\hat{g}(x,y)|^2 
+  |Z\hat{g}(x+p,y)|^2 = \frac{1}{p}$ a.e..
\item[\rm (iii)] For all $j =0,\ldots,\tfrac{L}{2p}-1$ and $y \in \ZZ_L$, we 
have $\sum_{l =0}^{\frac{L}{p}-1}\hat{g}(y+lp)\hat{g}(y+lp + 2jp) = \tfrac{1}{p}
\delta_{j,0}$.
\item[\rm (iv)] $\Wil(g,\Lambda,\CC^L)$ is an orthonormal basis for $\CC^L$.
\end{enumerate}
\end{proposition}

\begin{proof}
The proof, while lengthy, is very similar to the proof of
Proposition~\ref{eq_rec_Z}. In fact, with obvious adaptations, such
as using the normalized Haar measure on $\ZZ_L$, i.e.,
$m(E) = \frac{1}{L}\sum_{x \in \ZZ_L} 1_{E}(x)$ for all $E \subseteq \ZZ_L$, 
and replacing Zak
transforms and Fourier transforms by their corresponding finite
counterparts, the proof carries over almost line by line. We therefore
leave this part to the reader.
\end{proof}

Now we will turn our attention to general time--frequency lattices.
Here the situation is slightly more involved compared to $\ltZ$.
 
Let $\Lambda$ be a lattice in $\ZZ_L \times \ZZ_L$ with canonical generator 
matrix $A$ given by
\begin{equation*}
A=\begin{bmatrix} \frac{L}{2p} & b \\ 0 & p \end{bmatrix}.
\end{equation*}
Then we choose $\alpha,\beta,m_0,n_0 \in \ZZ$ such that 
\begin{equation}\label{minimal} 
\alpha \tfrac{L}{2p} m_0 + \alpha b n_0 + \beta p n_0
\end{equation}
attains its minimal positive value. 
Assume that there exists a choice of $\alpha,\beta,m_0,n_0$ such that
$(\tfrac{L}{2p}m_0 + bn_0)(pn_0) < 0$ and 
$(\alpha \tfrac{L}{2p})(\alpha b+\beta p)>0$. In the other
cases we have to change the signs of the later defined $\gamma$
and $\delta$ accordingly. In the following we will restrict to
the case where $|\alpha|=1$. For the
remainder of this section let $\alpha,\beta,m_0,n_0$ be defined
in this way. Now we
regard $\alpha$ and $\beta$ as elements of $\ZZ_L$. For the
sake of brevity we set
\[ c := \gcd(\alpha \tfrac{L}{2p},\alpha b+\beta p), \quad
d:= {\rm lcm}(\alpha \tfrac{L}{2p},\alpha b+\beta p),\]
and 
\[  s :={\rm gcd}( \tfrac{L}{2p}m_0 + bn_0,pn_0), \quad
t:= {\rm lcm}( \tfrac{L}{2p}m_0 + bn_0,pn_0).\]
The minimality condition for (\ref{minimal}) shows that
\begin{equation}\label{cs} 
\alpha \tfrac{L}{2p} m_0 + \alpha b n_0 + \beta p n_0 = c = s.
\end{equation}
We further define $\gamma,\delta \in \ZZ_L$ by
\[ \gamma := \frac{t}{\frac{L}{2p}m_0 + bn_0}\quad
\mbox{and} \quad \delta := -\frac{t}{pn_0}\]
and $\sigma \in {\rm Hom}(\ZZ_L \times \ZZ_L)$ by
\[ \sigma = \begin{bmatrix} \alpha & \beta \\ \gamma & \delta \end{bmatrix}.\]

\begin{proposition}
\label{meta_CL}
Let $\Lambda$ be a lattice in $\ZZ_L \times \ZZ_L$ with canonical generator 
matrix $A$ given by
\begin{equation*}
A=\begin{bmatrix} \frac{L}{2p} & b \\ 0 & p \end{bmatrix}.
\end{equation*}
Let $\sigma$ be defined
as in the preceding paragraph, and let $U$ on the space of $L$--periodic functions
on $\ZZ$ be defined by
\[ Uf(k) = \sum_{l \in \ZZ_L} 
f(\alpha k+\beta l) e^{-\pi i (\alpha \gamma k^2+ \beta \delta l^2)(L+1)/L}
e^{-2\pi i \beta \gamma kl/L}.\]
Then 
\[ g_{m\frac{L}{2p}+nb,np}(l) =
C(m,n)U(U^{-1}g)_{\sigma(m\frac{L}{2p}+nb,np)}(l),\]
where $C(m,n)=e^{-\pi i (\alpha \gamma m^2+ \beta \delta n^2)(L+1)/L}
e^{-2\pi i \beta \gamma mn/L}$.
\end{proposition}

Before moving on to the proof of this statement we point out that the 
operator $U$ in Proposition~\ref{meta_CL} is no longer a simple chirp 
operator as for the case $\ltZ$, cf.~Proposition~\ref{meta_Z}. This difference
and the different form of $\sigma$ necessitates a somewhat different
proof for the case $\CC^L$.

\begin{proof}
In order to apply Theorem \ref{Kai_Theorem} we need to check whether
$\sigma=(\sigma_1,\sigma_2)$ is symplectic. For this, we have to show that, for all 
$(x,y),(x',y') \in \ZZ_L \times \ZZ_L$, 
\begin{equation}\label{def_sym} e^{-2 \pi i \sigma_2(x,y) \sigma_1(x',y')/L}
e^{2 \pi i \sigma_2(x',y') \sigma_1(x,y)/L} = e^{-2 \pi i x' y/L}
e^{2 \pi i x y'/L}.\end{equation}
We have
\begin{eqnarray*}
\lefteqn{e^{-2 \pi i \sigma_2(x,y) \sigma_1(x',y')/L}
e^{2 \pi i \sigma_2(x',y') \sigma_1(x,y)/L}}\\
& = & e^{-2 \pi i (\gamma x + \delta y) (\alpha x' +\beta y') /L}
e^{2 \pi i (\gamma x' + \delta y') (\alpha x +\beta y)/L}\\
& = & e^{2 \pi i (\alpha\delta-\beta\gamma)(xy'-x'y)/L}
\end{eqnarray*}
and
\begin{eqnarray*} 
\alpha\delta-\beta\gamma
& = & -\frac{\alpha t}{pn_0} - \frac{\beta t}{\frac{L}{2p}m_0 + bn_0}\\
& = & \frac{-t}{(\frac{L}{2p}m_0 + bn_0)(pn_0)} (\alpha(\tfrac{L}{2p}m_0 + bn_0)
+ \beta p n_0).
\end{eqnarray*}
By (\ref{cs}),
\[ \alpha\delta-\beta\gamma = \frac{-st}{(\frac{L}{2p}m_0 + bn_0)(pn_0)}
= 1,\]
since $(\tfrac{L}{2p}m_0 + bn_0)(pn_0) < 0$. This proves that (\ref{def_sym}) is satisfied, which
shows that $\sigma$ is indeed symplectic. 
Moreover, a short computation analogous to the one in the proof of Proposition
\ref{meta_Z} shows that
\[ \zeta(k,l) = (\sigma^* \kappa_0 \sigma-\kappa_0)(k,l)
= (\alpha \gamma k + \beta \gamma l , \beta \gamma k + \beta \delta l).\]
Now it is easy to check (compare also \cite[Example 1.1.34 (iii)]{Kai99}) that
\[\psi : \ZZ_L^2 \to \TT, \quad \psi(k,l)
=e^{\pi i (\alpha \gamma k^2+ \beta \delta l^2)(L+1)/L}
e^{2\pi i \beta \gamma kl/L}\]
is a second degree character associated to $\zeta$. 
Applying Theorem \ref{Kai_Theorem} now finishes the proof.
\end{proof}

As in the discrete case we need to define a special bijective
map in order to give the definition of a Wilson basis associated 
with a lattice with arbitrary canonical generator matrix. 

\begin{lemma}
\label{def_varphiCL}
Let $\frac{L}{2p}, b \in \ZZ_L$ with $0 \le b < \frac{L}{2p}$ and
let $\alpha,\beta,m_0,n_0,c,d$ be defined as before.
Then the mapping $\varphi : \ZZ^2 \to \ZZ^2$ defined by
\[ \phi(m,n) = \left\{\begin{array}{ccl}
(m,n) & : & b=0,\\
(mm_0-\frac{2pd}{\alpha L}n,mn_0+\frac{d}{\alpha b + \beta p}n) & : & b \neq 0
\end{array}\right.\]
is bijective and we have
\[
\{(m \hspace*{-0.2cm} \mod \tfrac{L}{c},n \hspace*{-0.2cm} \mod 2c) : (m,n) \in \varphi^{-1}(\{0,\ldots,2p-1\}
\times \{0,\ldots,\tfrac{L}{p}-1\})\}\] 
\[= \{0,\ldots,\tfrac{L}{c}-1\} \times \{0,\ldots,2c-1\}.
\]
\end{lemma}

\begin{proof}
The proof of this lemma is very similar to the proof of
Lemma~\ref{def_varphi},
we therefore omit it.
\end{proof}

Note that the following
definition reduces to Definition \ref{defiwil_rectangleCL} in
the case of a diagonal canonical generator matrix. 

\begin{definition}
Let $\Lambda$ be a lattice in $\ZZ_L \times \ZZ_L$ with canonical generator 
matrix $A$ given by
\begin{equation*}
A=\begin{bmatrix} \frac{L}{2p} & b \\ 0 & p \end{bmatrix}.
\end{equation*}
Let $g$ be some $L$--periodic function on $\ZZ$, and let
$\varphi$ be defined as in Lemma \ref{def_varphiCL}.
Then the Wilson system 
$\Wil(g,\Lambda,\CC^L) = \{\psi_{m,n}\}_{(m,n) \in I}$, where $I = \{0,\ldots,p-1\}
\times \{0,\frac{L}{2p}\} \cup \{0,\ldots,2p-1\} \times \{1,\ldots,\frac{L}{2p}-1\}$, 
is given by
\[ \hspace*{1.8cm} \psi_{m,n}^\Lambda = 
g_{\varphi_1(2m,n)\frac{L}{2p}+\varphi_2(2m,n)b,\varphi_2(2m,n)p}, 
\hspace*{0.9cm} \mbox{if } m=0,\ldots,p-1,\;n=0,\tfrac{L}{2p},\]
and for $m=0,\ldots,2p-1$, $n=1,\ldots,\frac{L}{2p}-1,$
\begin{eqnarray*}
\psi_{m,n}^\Lambda & = & \tfrac{1}{\sqrt{2}} 
(g_{\varphi_1(m,n)\frac{L}{2p}+\varphi_2(m,n)b,\varphi_2(m,n)p}\\
&&
+ g_{\varphi_1(m,-n)\frac{L}{2p}+\varphi_2(m,-n)b,\varphi_2(m,-n)p}), \quad
\mbox{if } m+n \mbox{ even},\\
\psi_{m,n}^\Lambda & = & \tfrac{i}{\sqrt{2}} 
(g_{\varphi_1(m,n)\frac{L}{2p}+\varphi_2(m,n)b,\varphi_2(m,n)p}\\
&&
- g_{\varphi_1(m,-n)\frac{L}{2p}+\varphi_2(m,-n)b,\varphi_2(m,-n)p}), \quad
\mbox{if } m+n \mbox{ odd}.
\end{eqnarray*}
\end{definition}

The following theorem is the analog to Theorem~\ref{theo_general_Z} for
the space $\CC^L$.

\begin{theorem}
Let $\Lambda$ be a lattice in $\ZZ_L \times \ZZ_L$ with canonical generator 
matrix $A$ given by
\begin{equation*}
A=\begin{bmatrix} \frac{L}{2p} & b \\ 0 & p \end{bmatrix}.
\end{equation*}
Let $g$ be some $L$--periodic function on $\ZZ$ such that $\widehat{U^{-1}g}$ is real--valued, 
let $M:=2p$, $N:=\frac{L}{p}$, $q := \frac{L}{2c}$,
$\tilde{M}:=2q$, $\tilde{N}:=\frac{L}{q}$, and let $U$ and $\varphi$ be defined as in 
Proposition \ref{meta_CL} and Lemma \ref{def_varphiCL}, respectively.
Then the following conditions are equivalent.
\begin{enumerate}
\item[\rm (i)] $\{g_{m\frac{L}{2p}+nb,np}\}_{m=0,\ldots,M-1,n=0,\ldots,N-1}$ is a tight 
frame for $\CC^L$ with frame bound $2$.
\item[\rm (ii)] $\{(U^{-1}g)_{m\frac{L}{2q},nq}\}_{m=0,\ldots,\tilde{M}-1,
n=0,\ldots,\tilde{N}-1}$ 
is a tight frame for $\CC^L$ with frame bound $2$.
\item[\rm (iii)] $\Wil(U^{-1}g,\frac{L}{2q}\ZZ_L  \times q\ZZ_L, \CC^L)$ is an orthonormal basis for $\CC^L$.
\item[\rm (iv)] $\Wil(g,\Lambda,\CC^L)$ is an orthonormal basis for $\CC^L$.
\end{enumerate}
\end{theorem}

\begin{proof}
Let $\sigma$ be defined as before and
let $\varphi = (\varphi_1,\varphi_2)$.
Then we compute
\begin{eqnarray*}
\lefteqn{\sigma(\varphi_1(m,n)\tfrac{L}{2p}+\varphi_2(m,n)b,\varphi_2(m,n)p)}\\
& = & (\alpha(\varphi_1(m,n)\tfrac{L}{2p}+\varphi_2(m,n)b)
+ \beta \varphi_2(m,n)p,
\gamma(\varphi_1(m,n)\tfrac{L}{2p}+\varphi_2(m,n)b) + \delta \varphi_2(m,n)p)\\
& = & (m(\alpha \tfrac{L}{2p} m_0 + \alpha b n_0 + \beta p n_0)
+n(-\alpha \tfrac{L}{2p}\tfrac{2pd}{\alpha L}+\alpha b\tfrac{d}{\alpha b + \beta p}
+\beta p\tfrac{d}{\alpha b + \beta p}),\\
& & m(\gamma \tfrac{L}{2p} m_0 + \gamma b n_0 + \delta p n_0)
+n(-\gamma \tfrac{L}{2p}\tfrac{2pd}{\alpha L}+\gamma b\tfrac{d}{\alpha b + \beta p}
+\delta p\tfrac{d}{\alpha b + \beta p}))\\
& = & (mc,m(\tfrac{t}{\frac{L}{2p}m_0 + bn_0}(\tfrac{L}{2p} m_0 +b n_0)-
\tfrac{t}{pn_0} p n_0)\\
&&
+n(-\tfrac{t}{\frac{L}{2p}m_0 + bn_0}\tfrac{d}{\alpha}
+\tfrac{t}{\frac{L}{2p}m_0 + bn_0}\tfrac{bd}{\alpha b + \beta p}
-\tfrac{t}{pn_0}\tfrac{pd}{\alpha b + \beta p}))\\
& = & (mc,ntd(\tfrac{-\alpha b p n_0 -\beta p^2 n_0 + \alpha b p n_0
-\alpha p \frac{L}{2p}m_0 - \alpha pbn_0}{\alpha (\frac{L}{2p}m_0 + bn_0)(\alpha b + \beta p)pn_0}))\\
& = & (mc,-ntdp(\tfrac{\alpha \frac{L}{2p}m_0 +\alpha bn_0 
+ \beta p n_0}{\alpha (\frac{L}{2p}m_0 + bn_0)(\alpha b + \beta p)pn_0}))\\
& = & (mc,-ntp\tfrac{cd}{\alpha (\frac{L}{2p}m_0 + bn_0)(\alpha b + \beta p)pn_0}))\\
& = & (mc,-nt\tfrac{\frac{L}{2}}{(\frac{L}{2p}m_0 + bn_0)pn_0})\\
& = & (mc,n\tfrac{L}{2c})),
\end{eqnarray*}
where in the last step we used (\ref{cs}) and $st=-(\frac{L}{2p}m_0 + bn_0)pn_0$.
Since $|\alpha|=1$, we have $c=s=\gcd(\alpha \tfrac{L}{2p},\alpha b+\beta p)$
is a factor of $\tfrac{L}{2p}$ and hence of $\tfrac{L}{2}$.
Using Lemma \ref{def_varphiCL}, the equivalence of 
(i) and (ii) follows immediately from Proposition \ref{meta_CL},
since $U$ is unitary and $|C(m,n)|=1$.
Proposition \ref{eq_rec_CL} proves (ii) $\Leftrightarrow$ (iii). 
Therefore it remains to prove the equivalence of (iii) and (iv). 
For this, we will use the following 
implication of Proposition \ref{meta_CL}:
\begin{eqnarray*} 
U(U^{-1}g)_{m\frac{L}{2q},nq}
&=& U(U^{-1}g)_{\sigma(\varphi_1(m,n)\frac{L}{2p}+\varphi_2(m,n)b,\varphi_2(m,n)p)}\\
&=& C(\varphi(m,n))^{-1}g_{\varphi_1(m,n)\frac{L}{2p}+\varphi_2(m,n)b,\varphi_2(m,n)p}.
\end{eqnarray*}
An easy but tedious calculation shows that $C(\varphi(m,n))^{-1}$ does not depend on
the sign of $n$. 
Using now the definition of a Wilson basis, the fact that $U$ is a unitary
operator, and the fact that $|C(\varphi(m,n))|=1$ yields the result. 
\end{proof}

Tight Gabor frames in $\Cst^L$ can be constructed in the same
way as for $\ltZ$ and $\LtR$ by using the ``inverse square root trick''.
Furthermore, it has been shown in~\cite{Kai04} that for properly localized 
windows the dual window constructed in $\Cst^L$ by ``sampling and 
periodization'' of the frame $\{g_{ma,nb}\}$ converges to the dual 
window $S^{-1}g$ with increasing sampling rate and increasing periodization 
interval, see~\cite{Kai04} for details. This result can be easily extended
to tight windows. We refer also to~\cite{Jan97a,Str01,CS03} for related 
results and leave the details to the reader. To obtain tight Gabor frames 
in $\Cst^L$ that satisfy the required conditions of the theorem above and have 
good time-frequency localization one can thus essentially proceed analogous
to the example at the end of Section~\ref{s:cont}.

\section{Conclusion} \label{s:conc}

We have demonstrated that orthonormal Wilson bases for $\LtR$ (with excellent 
time-frequency localization) can be constructed for general time-fre\-quen\-cy
lattices. Of course any numerical implementation has to be done in a discrete 
setting. Somewhat longer proofs establish a similar result for the 
spaces $\ltZ$ and $\Cst^L$ for non-rectangular time-frequency lattices.
The approach based on metaplectic transforms used in this paper suggests that
the main results can be extended to the setting of symplectic
time-frequency lattices on general locally compact abelian groups.

Furthermore, our results imply that from a practical viewpoint 
it is indeed possible to extend OQAM-OFDM or cosine-modulated filter banks 
to general time-frequency lattices.  Moreover,
we expect that the benefits of using general time-frequency lattices will 
be even more pronounced for images and higher-dimensional signals. 
Our expectation is based on the fact that in the theory of sphere packings
(and sphere coverings) the advantages of the optimal sphere packing over
the packing associated with the rectangular lattices increases 
significantly with the dimension of the space~\cite{CS93b}.

An interesting research problem is thus to investigate how to extend
the results in this paper to $L^2(\RR^d)$ for non-symplectic
lattices as well as to find optimal time-frequency 
lattices in $\RR^{2d}$ for $d>1$. One possibility to define an ``optimal''
time-frequency lattice is to fix the function $g$ to be a Gaussian, say,
and then find that time-frequency lattice of fixed density which minimizes
the condition number of the associated Gabor frame operator as
indicated in~\cite{SB03}.

\section*{Acknowledgments}

We thank the referees for valuable comments and suggestions which lead
to an improvement of the results and presentation in this paper.


\begin{thebibliography}{10} 
\bibitem{AT90}
{\sc L.~Auslander and R.~Tolimieri},
{\em On finite {G}abor expansions of signals},
in Signal Processing, Part I: Signal Processing Theory, 
L.~Auslander, T.~Kailath, and S.K. Mitter, eds., Springer Verlag,
  New York, 1990, pp.~13--23, IMA vol. 22, lectures from IMA Program, summer 1988.

\bibitem{BCM03}
{\sc J.J. Benedetto, W.~Czaja, and A.Y. Maltsev},
{\em The {B}alian-{L}ow theorem for the symplectic form on {$\mathbb{R}^{2d}$}},
J. Math. Phys., 44 (2003), pp.~1735--1750.

\bibitem{BHW95}
{\sc J.J. Benedetto, C.~Heil, and D.F. Walnut},
{\em Differentiation and the {B}alian--{L}ow theorem},
J.\ Four.\ Anal.\ Appl., 1 (1995), pp.~355--403.

\bibitem{Bit03}
{\sc K. Bittner},
{\em Wilson bases on the interval},
in Advances in Gabor Analysis, H.G. Feichtinger and T.~Strohmer, eds., 
Appl. Numer. Harmon. Anal.,  Birkh\"auser,
Boston, MA, 2003, pp.~197--221.

\bibitem{Bol02}
{\sc H.~B\"olcskei},
{\em Orthogonal frequency division multiplexing based on offset {QAM}},
in Advances in Gabor Analysis, H.G. Feichtinger and T.~Strohmer, eds., 
Birkh\"auser, Boston, MA, 2002, pp.~321--352. 

\bibitem{BFG97}
{\sc H.~B{\"o}lcskei, H.G.~Feichtinger, K.~Gr{\"o}chenig, and F.~Hlawatsch.}
\newblock Discrete-time Wilson expansions.
\newblock {\em Proc. IEEE Int. Sympos. Time-Frequency Time-Scale Analysis}, 
Paris (France), pp. 525-528, June 1996.

\bibitem{BH98}
{\sc H.~B{\"o}lcskei and F.~Hlawatsch},
{\em Oversampled modulated filter banks},
in Gabor {A}nalysis and {A}lgorithms: {T}heory and {A}pplications, H.G. Feichtinger and T.~Strohmer, eds., 
Birkh\"auser, Boston, MA, 1998, pp.~295--322.

\bibitem{CS93b}
{\sc J.H. Conway and N.J.A. Sloane}, editors.
\newblock {\em Sphere Packings, Lattices and Groups}.
\newblock Grundlehren der mathematischen Wissenschaften. Springer Verlag,
New York, Berlin, Heidelberg, 1993.


\bibitem{Chr03}
{\sc O.~Christensen},
{\em An introduction to frames and Riesz bases}, Birkh\"auser, Boston, 2003.

\bibitem{CS03}
O.~Christensen and T.~Strohmer.
\newblock Methods for approximation of the inverse ({G}abor) frame
operator.
\newblock In {\em Advances in Gabor analysis}, pages 171--195. Birkh\"auser
  Boston, Boston, MA, 2003.
                                                                                

\bibitem{DJJ91}
{\sc I.~Daubechies, S.~Jaffard, and J.L. Journ\'{e}},
{\em A simple {W}ilson orthonormal basis with exponential decay},
SIAM J. Math. Anal., 22 (1991), pp.~554--572.

\bibitem{FCS95}
{\sc H.G. Feichtinger, O.~Christensen, and T.~Strohmer},
{\em A group-theoretical approach to {G}abor analysis},
Optical Engineering, 34 (1995), pp.~1697--1704.

\bibitem{FK98}
{\sc H.G. Feichtinger and W.~Kozek},
{\em Quantization of {TF}--lattice invariant operators on elementary {LCA}
  groups},
in Gabor {A}nalysis and {A}lgorithms: {T}heory and {A}pplications, H.G. Feichtinger and T.~Strohmer, eds., 
Birkh\"auser, Boston, MA, 1998, pp.~233--266.

\bibitem{FS98}
{\sc H.G. Feichtinger and T.~Strohmer}, eds.,
{\em {G}abor Analysis and Algorithms: Theory and Applications},
Birkh\"auser, Boston, 1998.

\bibitem{FS03}
{\sc H.G. Feichtinger and T.~Strohmer}, eds.,
{\em Advances in {G}abor analysis},
Birkh\"auser Boston Inc., Boston, MA, 2003.

\bibitem{FAB95}
{\sc B.~Le Floch, M.~Alard, and C.~Berrou},
{\em Coded orthogonal frequency division multiplex},
Proc. of IEEE, 83 (1995), pp.~982--996.

\bibitem{Fol89}
{\sc G.B. Folland},
{\em Harmonic Analysis in Phase Space},
Annals of Math, Studies. Princeton Univ. Press, Princeton (NJ), 1989.

\bibitem{fol} 
{\sc G.B. Folland}, 
{\em A course in abstract harmonic analysis},
CRC Press, Baco Raton, 1995.

\bibitem{Gro98}
{\sc K.~Gr{\"o}chenig},
{\em Aspects of {G}abor analysis on locally compact abelian groups},
in Gabor {A}nalysis and {A}lgorithms: {T}heory and {A}pplications, H.G. Feichtinger and T.~Strohmer, eds., 
Birkh\"auser, Boston, MA, 1998, pp.~211--231.

\bibitem{Gro01}
{\sc K.~{Gr\"ochenig}},
{\em Foundations of Time-Frequency Analysis},
Birkh\"auser, Boston, 2001.

\bibitem{GHH02}
{\sc K.~Gr{\"o}chenig, D.~Han, C.~Heil, and G.~Kutyniok},
{\em The {B}alian--{L}ow theorem for symplectic lattices in higher
  dimensions},
Appl.\ Comp.\ Harm.\ Anal., 13 (2002), pp.~169--176.

\bibitem{Her1850}
{\sc C. Hermite}, 
{\em Extraits de lettres de M. Ch. Hermite \`{a} M. Jacobi
sur differents objets de la th\'{e}orie des nombres, Deuxi\`{e}me lettre},
Reine Angewandte Mathematik, 40 (1850), pp.~279--290.

\bibitem{hew} 
{\sc E. Hewitt and K.A. Ross}, {\em Abstract harmonic analysis I, II},
Springer-Verlag, Berlin/Hei\-del\-berg/New York, 1963/1970.

\bibitem{Hua82}
{\sc L.K. Hua},
{\em Introduction to number theory},
Springer-Verlag, Berlin, 1982.

\bibitem{Jan97a}
A.J.E.M. Janssen.
\newblock From continuous to discrete {W}eyl-{H}eisenberg frames through
  sampling.
\newblock {\em J. Fourier Anal. Appl.}, 3(5):583--596, 1997.
                                                                                


\bibitem{Kai99} 
{\sc N. Kaiblinger}, {\em Metaplectic representation, eigenfunctions
of phase space shifts, and Gelfand--Shilov spaces for lca groups}, 
Ph.D. thesis, University of Vienna, 1999.

\bibitem{Kai04} 
{\sc N. Kaiblinger}, 
{\em Approximation of the Fourier transform and the dual Gabor window}, 
J.~Fourier Anal. Appl., to appear.

\bibitem{KK} {\sc E. Kaniuth and G. Kutyniok}, 
{\em Zeros of the Zak transform on locally compact abelian groups}, 
Proc. Amer. Math. Soc. {\bf 126} (1998), 3561--3569. 

\bibitem{KM98}
{\sc W.~Kozek and A.~Molisch},
{\em Nonorthogonal pulseshapes for multicarrier communications in doubly
  dispersive channels},
IEEE J. Sel. Areas Comm., 16 (1998), pp.~1579--1589.

\bibitem{Kut00} {\sc G. Kutyniok}, {\em Time-frequency analysis on locally compact
groups}, Ph.D. thesis, University of Paderborn, 2000.

\bibitem{Str97a}
{\sc T.~Strohmer},
{\em Numerical algorithms for discrete {G}abor expansions},
in Gabor {A}nalysis and {A}lgorithms: {T}heory and {A}pplications, H.G. Feichtinger and T.~Strohmer, eds., 
Birkh\"auser, Boston, MA, 1998, pp.~267--294.

\bibitem{Str01}
{\sc T.~Strohmer},
{\em Approximation of dual {G}abor frames, window decay, and wireless
  communications},
Appl.\ Comp.\ Harm.\ Anal., 11 (2001), pp.~243--262.

\bibitem{SB03}
{\sc T.~Strohmer and S.~Beaver},
{\em Optimal {OFDM} system design for time-frequency dispersive channels},
IEEE Trans.\ Comm., 51:7 (2003), pp.~1111--1122.

\bibitem{TA98}
{\sc R.~Tolimieri and M.~An},
{\em Time-frequency representations},
Applied and Numerical Harmonic Analysis, Birkh\"auser Boston Inc.,
  Boston, MA, 1998.

\bibitem{Wil87}
{\sc K.G. Wilson},
{\em Generalized {W}annier functions}, (1987), preprint.

\end{thebibliography}
\end{document}